\documentclass[10pt]{article}
\usepackage{amsfonts}
\usepackage{amsmath}
\usepackage{mathrsfs}
\usepackage{float}
\usepackage{graphicx}
\usepackage{subfigure}
\usepackage{color}
\usepackage{mathrsfs,amscd,amssymb,amsthm,amsmath,bm,graphicx,psfrag,subfigure,url}
\input{epsf}
\makeatletter

\renewcommand{\@seccntformat}[1]{{\csname the#1\endcsname}{\normalsize .}\hspace{.5em}}
\makeatother

\usepackage{indentfirst}

\def \[{\begin{equation}}
\def \]{\end{equation}}

\newtheorem{thm}{Theorem}[section]

\newtheorem{defi}{Definition}

\newtheorem{lem}[thm]{Lemma}
\newtheorem{cor}[thm]{Corollary}

\newenvironment{wst}
{\setlength{\leftmargini}{1.5\parindent}
 \begin{itemize}
 \setlength{\itemsep}{-1.1mm}}
{\end{itemize}}

\begin{document}

\setlength{\baselineskip}{17pt}
\begin{center}{\Large \bf On the resistance distance and Kirchhoff index of a linear
hexagonal (cylinder) chain\footnote{Financially supported by the National Natural Science Foundation of China (Grant Nos. 11671164, 11271149).}}
\vspace{4mm}

{\large Sumin Huang,\ Shuchao Li $^{a,}$\footnote{Corresponding author.\\
\hspace*{5mm}E-mail: sumin2019@sina.com (S. Huang),
\ lscmath@mail.ccnu.edu.cn (S.C.
Li)}}\vspace{2mm}

Faculty of Mathematics and Statistics,  Central China Normal
University, Wuhan 430079, PR China
\end{center}

\noindent {\bf Abstract}:\ The resistance between two nodes in some resistor networks has been studied extensively by mathematicians and physicists. Let $L_n$ be a linear hexagonal chain with $n$\, 6-cycles. Then identifying the opposite lateral edges of $L_n$ in ordered way yields the linear hexagonal cylinder chain, written as $R_n$. We obtain explicit formulae for the resistance distance $r_{L_n}(i, j)$ (resp. $r_{R_n}(i,j)$) between any two vertices $i$ and $j$ of $L_n$ (resp. $R_n$). To the best of our knowledge $\{L_n\}_{n=1}^{\infty}$ and $\{R_n\}_{n=1}^{\infty}$ are two nontrivial families with diameter going to $\infty$ for which all resistance distances have been explicitly calculated. We determine the maximum and the minimum resistance distances in $L_n$ (resp. $R_n$). The monotonicity and some asymptotic properties of resistance distances in $L_n$ and $R_n$ are given. As well we give formulae for the Kirchhoff indices of $L_n$ and $R_n$ respectively.

\vspace{2mm} \noindent{\bf Keywords}: Resistance distance; Kirchhoff index; Moore-Penrose inverse

\vspace{2mm}\noindent{2010 Mathematics Subject Classification. 05C50.}

 {\setcounter{section}{0}}
\section{\normalsize Introduction}\setcounter{equation}{0}
The resistance distance (known also as the effective resistance) of a graph is one important measure of quantifying structural
properties for the given graph. The resistance is not only suggested \cite{13,14} to be a central concept in electronic circuit theory, but also has widespread utility in physics, engineering, mathematics, chemistry and computer sciences. It has been shown \cite{15,16,17} that the escape probability, the first passage time, the cover cost and the commute time of random walks have closely relation with the resistance. For more advances one may be referred to \cite{Georgakopoulos1,Georgakopoulos2,Huang1,Huang2} and the references cited in.

The computation of two-point resistance of a graph is a classical problem in electric circuit theory, which attracts much attention \cite{12}.
Gervacio \cite{24} obtained an explicit expression for the resistance between any pair of vertices in the complete $n$-partite graph.  Based on the Gervacio��s method, Jiang and Yan \cite{25} obtained the closed formula of the resistance in so-called ring network graphs. Cinkir \cite{26} obtained explicit formulae for Kirchhoff index and resistances between vertices of linear polyomino chain. Shi and Chen \cite{22} used a new method to obtained explicit formulae for resistances between vertices of linear polyomino chain, and determine the largest and the smallest resistances in linear polyomino chain. For the wheel and the fan, resistance distance between any two vertices has been calculated explicitly as a function of the number of vertices in the graph \cite{Bapat,Yang1}. Vaskouski and Zadorozhnyuk \cite{Vaskouski} studied the resistance distance between any two vertices in Cayley graphs on symmetric groups.

Recently, concerns were raised that resistance distance fails a number of desirable properties of a distance function for
certain random geometric graphs \cite{Ulrike}. For these graphs they obtain the asymptotic result that
$$
 r_G(i,j) \approx \frac{1}{{\rm deg}(i)}+\frac{1}{{\rm deg}(j)}.
$$
Note that the value of $r_G(i, j)$ here depends only on the degrees of vertices $i$ and $j$, they concluded that $r_G(i, j)$ is completely
meaningless as a distance function on these large geometric graphs.

Clearly, the preceding result does not hold for some classes of graphs. For trees $r_G(i, j) = d_G(i, j)$, so $r_G(i, j)$ is still a
distance function. We know that resistance distance has been calculated for a number of special graphs as above. However, there are no widely available effective computational tools to compute the effective resistance of a graph of reasonable size. It still seems to be a paucity of results for infinite classes of graphs. Barrett, Evans and Francis \cite{barrett} investigate the resistance distance for an infinite class of 2-trees.

Motivated by \cite{barrett,26,22}, we study other two infinite classes of graphs, i.e., the linear hexagonal chain and the linear hexagonal cylinder chain, for which the effective resistance retains all desirable properties of a distance function.

This paper is organized by the following way. In Section 2, we give some necessary definitions and preliminary results. In Section 3,  we first obtain explicit formulae for the resistance distance between any two vertices in the linear hexagonal chain. Then we determine the largest and the smallest resistances in the linear hexagonal chain. The monotonicity and asymptotic property of resistances in $L_n$ are discussed. In Section 4, we first obtain explicit formulae for the resistance distance between any two vertices in the linear hexagonal cylinder chain. Then we determine the largest and the smallest resistances in the linear hexagonal cylinder chain. As well the monotonicity and asymptotic property of resistances in $R_n$ are discussed. Based on our obtained results in this paper we obtained the formulae for the Kirchhoff indices (i.e., $K\!f(G)=\sum_{\{u,v\}\subseteq V_G}r_G(u,v)$) of $L_n$ and $R_n$ in the last section. It is interesting to see that $\frac{K\!f(R_n)}{K\!f(L_n)}\rightarrow \frac{1}{2}$ as $n\rightarrow \infty.$

\section{\normalsize  Some definitions and preliminary results}\setcounter{equation}{0}

In this section, we give some necessary definitions and preliminary results. A graph is denoted by $G=(V_G, E_G)$, where $V_G$ is the vertex set and $E_G$ is the edge set. The \textit{order} of $G$ is the number $|V_G|$ of its vertices, and the \textit{size} is the number $m=|E_G|$ of its edges.

There are, however, chemically interesting unbranched polycyclic polymers which are uniform. This means that they are composed of cycles of uniform lengths.
Probably the best known and the most relevant are the linear hexagonal chain and the linear hexagonal cylinder chain (or hexagonal cylinder chain, for short), consisting of $n$\, $6$-cycles, which are depicted in Fig.~1.
\begin{figure}[!ht]
\centering
\psfrag{r}{$i$}
\psfrag{2}{$\cdots$}
\psfrag{k}{$n$}
\psfrag{f}{$L_n$}
\psfrag{t}{$R_n$}
\psfrag{h}{$R_6$}
\psfrag{7}{$R_1$}
\psfrag{8}{$R_2$}
\psfrag{u}{$1$}
\psfrag{x}{$p_{i-1}$}
\psfrag{y}{$u_{i-1}$}
\psfrag{z}{$p_i$}
\psfrag{5}{$q_{i-1}$}
\psfrag{6}{$v_{i-1}$}
\psfrag{7}{$q_i$}
  \includegraphics[width=120mm]{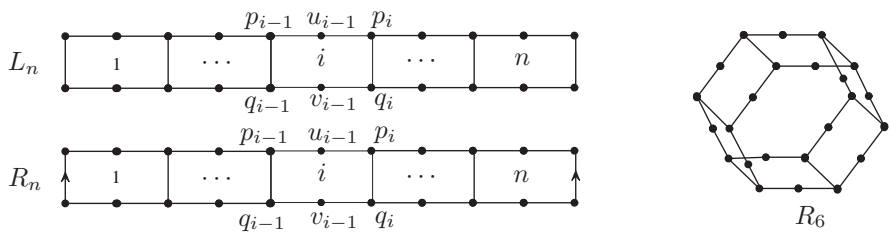}
  \caption{The linear hexagonal chain $L_n,$ and the linear hexagonal cylinder chains $R_n$ and $R_6$.}
\end{figure}

There are many techniques which are employed to calculate resistance distance, including the well-known series and parallel rules and the $\Delta$-$Y$ transformation, which are list in what follows.
\begin{defi}[Series Transformation] Let $N_1, N_2$, and $N_3$ be nodes in a graph where $N_2$ is adjacent to only $N_1$ and $N_3$. Moreover,
let $R_a$ equal the resistance between $N_1$ and $N_2$ and $R_b$ equal the resistance between node $N_2$ and $N_3$. A series transformation
transforms this graph by deleting $N_2$ and setting the resistance between $N_1$ and $N_3$ equal to $R_a + R_b$.
\end{defi}
\begin{defi}[Parallel Transformation] Let $N_1$ and $N_2$ be nodes in a multi-edged graph where $e_1$ and $e_2$ are two edges
between $N_1$ and $N_2$ with resistances $R_a$ and $R_b$, respectively. A parallel transformation transforms the graph by deleting
edges $e_1$ and $e_2$ and adding a new edge between $N_1$ and $N_2$ with edge resistance $r=(\frac{1}{R_a}+\frac{1}{R_b})^{-1}$
\end{defi}
A $\Delta$-$Y$ transformation is a mathematical technique to convert resistors in a triangle formation to an equivalent system
of three resistors in a $``Y"$ format as illustrated in Fig. 2. We formalize this transformation below.
\begin{figure}[!ht]
\centering
\psfrag{1}{$u$}
\psfrag{2}{$v$}
\psfrag{3}{$w$}
\psfrag{4}{$R_c$}
\psfrag{5}{$R_a$}
\psfrag{6}{$R_b$}
\psfrag{7}{$R_1$}
\psfrag{8}{$R_2$}
\psfrag{9}{$R_3$}
  \includegraphics[width=90mm]{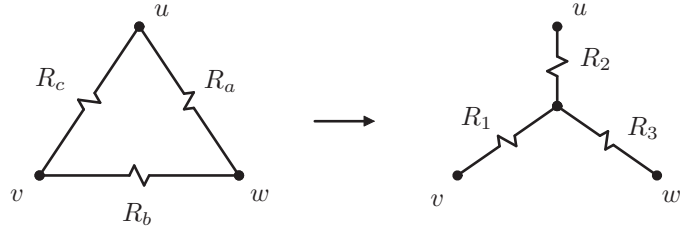}
  \caption{$\Delta$ and $Y$ circuits with vertices labeled as in Definition \ref{defi3}.}
  \label{fig10}
\end{figure}

\begin{defi}[$\Delta$-$Y$ Transformation]\label{defi3}
Let $N_1,N_2,N_3$ be nodes and $R_a,R_b$ and $R_c$ be given resistances as shown in Fig. $2$. The
transformed circuit in the $``Y"$ format as shown in Fig. $2$ has the following resistances:
\begin{eqnarray*}
  R_1=\frac{R_bR_c}{R_a+R_b+R_c}, \ \  R_2=\frac{R_aR_c}{R_a+R_b+R_c}, \ \
  R_3=\frac{R_aR_b}{R_a+R_b+R_c}.
\end{eqnarray*}
\end{defi}
\begin{lem}[\cite{6}]\label{lem1}
   Series transformations, parallel transformations, and $\Delta$-$Y$ transformations yield equivalent circuits.
\end{lem}

Further on we need the the following lemma. 
\begin{lem}[\cite{Bapat2,Klein}]\label{lem2}
Assume that $G$ is a graph with discrete Laplacian $L$. Use $L^{+}$ to denote the pseudo inverse of $L$, then we have
$$r(p,q)=l^{+}_{pp}+l^{+}_{qq}-2 l^{+}_{pq},$$
for every vertices p and q of G, where $l^{+}_{pp}$, $l^{+}_{qq}$ and $l^{+}_{pq}$ are the elements of the matrix $L^{+}$.
\end{lem}


\section{\normalsize The effective resistance in $L_n$}
\subsection{\normalsize Determining the effective resistance between any two vertices in $L_n$}
In this subsection, we determine the resistance distance for every pair of vertices of $L_n$. We label the $i$th 6-cycle of $L_n$ as $p_{i-1}, u_{i-1}, p_i, q_i, v_{i-1}, q_{i-1}$ for $i=1,2,\ldots, n$ (see Fig. 1). In our context, we abbreviate $r_{L_n}(u,v)$ to $r(u,v)$ for any pair of vertices $u,v$ in $L_n$. The order of $L_n$ is $4n+2$, whereas its size is $5n+1$.

First, we compute the effective resistance between $p_n$ and $q_n$. Let $z_n=r(p_n,q_n)$. We can express $z_{n+1}$ in terms of $z_n$ by the parallel circuit reduction in Fig. 3, that is:
\begin{align*}
 z_{n+1}=&\frac{1}{1+\frac{1}{z_n+4}} = \frac{z_n+4}{z_n+5},\ \ \ \  z_0=1.
\end{align*}
\begin{figure}[ht!]
\begin{center}
  \psfrag{a}{$p_n$}
  \psfrag{b}{$u_n$}
  \psfrag{c}{$p_{n+1}$}
  \psfrag{1}{$q_n$}
  \psfrag{2}{$v_n$}
  \psfrag{3}{$q_{n+1}$}
  \psfrag{z}{$z_n$}
  \psfrag{A}{$L_{n+1}$}
  \psfrag{B}{$\Rightarrow$}
  \includegraphics[width=60mm]{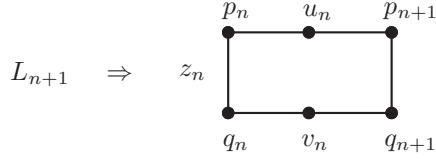}\\
  \caption{Graph $L_{n+1}$ and its simplified circuit.}
\end{center}
\end{figure}
Solving this recurrence relation gives
\[\label{eq:3.1}
  z_n=-(2+2\sqrt{2})+\frac{4\sqrt{2}}{1-(3-2\sqrt{2})^{2n+2}}.
\]

Set $\alpha:=3-2\sqrt{2}$. Then (\ref{eq:3.1}) can be rewritten as
$$z_n=-(2+2\sqrt{2})+\frac{4\sqrt{2}}{1-\alpha^{2n+2}}.$$

Our next aim is to find the effective resistances between any pair among $\{p_0, q_0, p_n\}$. Let $x_n=r(p_n,p_0),\, y_n=r(p_n,q_0)$.
By repeatedly using Lemma \ref{lem1}, we may obtain the simplified circuit of $L_n$ as depicted in Fig. 4. Hence, $x_{n-1}=A+C$, $y_{n-1}=B+C$, $z_{n-1}=A+B$. Then using parallel and series circuit reductions yields
\begin{align}
x_n=&\frac{1}{\frac{1}{A+2}+\frac{1}{B+3}}+C\notag\\
=&\frac{(x_{n-1}-y_{n-1}+z_{n-1}+4)(y_{n-1}-x_{n-1}+z_{n-1}+6)}{4(z_{n-1}+5)}+\frac{x_{n-1}+y_{n-1}-z_{n-1}}{2},\label{eq:3.2}\\
y_n=&\frac{1}{\frac{1}{B+2}+\frac{1}{A+3}}+C\notag\\
=& \frac{(x_{n-1}-y_{n-1}+z_{n-1}+6)(y_{n-1}-x_{n-1}+z_{n-1}+4)}{4(z_{n-1}+5)}+\frac{x_{n-1}+y_{n-1}-z_{n-1}}{2}.\label{eq:3.3}
\end{align}
This gives $x_n-y_n=\frac{x_{n-1}+y_{n-1}}{z_{n-1}+5}$. Let $t_n=x_n-y_n$. Then we have
$$\text{$t_n=\frac{t_{n-1}}{z_{n-1}+5},$\ if $n\geq 1$\  and \  $t_0=-1.$} $$

\begin{figure}[h!]
\begin{center}
  \psfrag{r}{$p_0$}\psfrag{2}{$q_0$}\psfrag{3}{$q_n$}\psfrag{4}{$C$} \psfrag{A}{$L_n$}\psfrag{B}{$\Rightarrow$}
   \psfrag{c}{$p_{n-1}$}\psfrag{1}{$q_{n-1}$}\psfrag{d}{$u_{n-1}$}\psfrag{2}{$v_{n-1}$} \psfrag{e}{$p_n$}\psfrag{c}{$p_{n-1}$}
   \psfrag{d}{$u_{n-1}$}\psfrag{s}{$B$}\psfrag{t}{$A$} \psfrag{H}{$L_n$}
  \includegraphics[width=80mm]{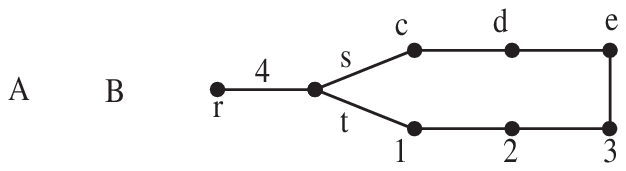}\\
  \caption{Graph $L_n$ and its simplified circuit.}
\end{center}
\end{figure}
Note that $z_n+5=\frac{(3+2\sqrt{2})^{n+2}-(3-2\sqrt{2})^{n+2}}{(3+2\sqrt{2})^{n+1}-(3-2\sqrt{2})^{n+1}}$. Hence,
\begin{align}\label{eq:3.4}
t_n=-\prod_{i=0}^{n-1}\frac{1}{z_i+5}=\frac{-4\sqrt{2}\alpha^{n+1}}{1-\alpha^{2n+2}}.
\end{align}
As $x_n=y_n+t_n$, together with (\ref{eq:3.2}), (\ref{eq:3.4}), and doing some algebra, (\ref{eq:3.3}) becomes
$$y_n=\frac{(t_{n-1}+z_{n-1}+6)(-t_{n-1}+z_{n-1}+4)}{4(z_{n-1}+5)}+\frac{t_{n-1}-z_{n-1}}{2}+y_{n-1}.$$
Note that $z_n-t_n=\frac{2\sqrt{2}-2+(2+2\sqrt{2})\alpha^{n+1}}{1-\alpha^{n+1}}=\frac{-t_{n-1}+z_{n-1}+4}{z_{n-1}+5}$. Hence,
\begin{align}
y_n=&\frac{z_n-t_n}{2} \frac{t_{n-1}+z_{n-1}+6}{2}-\frac{z_{n-1}-t_{n-1}}{2}+y_{n-1}\notag\\
=&\frac{z_n-t_n}{2}  \frac{t_{n-1}+z_{n-1}+4}{2}+\frac{z_n-t_n}{2}-\frac{z_{n-1}-t_{n-1}}{2}+y_{n-1}.\label{eq:3.5}
\end{align}
It is straightforward to check that $\frac{z_n-t_n}{2}  \frac{t_{n-1}+z_{n-1}+4}{2}=1$. Then (\ref{eq:3.5}) becomes
\begin{align*}
y_n=&\frac{z_n-t_n}{2}-\frac{z_{n-1}-t_{n-1}}{2}+y_{n-1}+1
=\frac{2\sqrt{2}}{1-\alpha^{n+1}}-1-\sqrt{2}+n.
\end{align*}
Therefore, we can compute $x_n$ by $y_n$ and $t_n$. That is,
$$x_n=t_n+y_n=\frac{2\sqrt{2}}{1+\alpha^{n+1}}-1-\sqrt{2}+n.$$

Next, we determine the formulae for $r(p_n,p_i)$, $r(p_n,q_i)$ and $r(p_i,q_i)$, where $n>i>0$. In fact, $L_n$ can be simplified to a $Y$-shaped graph as depicted in Fig. 5, where $M, N, K$ are the resistances in the $Y$-shaped graph. Thus, we have $M+N=x_{n-i-1}$, $M+K=y_{n-i-1}$, $N+K=z_{n-i-1}$. Equivalently,
\begin{align*}
M&=\frac{x_{n-i-1}+y_{n-i-1}-z_{n-i-1}}{2}=n-i-1,\ \ \ \ \
N=\frac{x_{n-i-1}-y_{n-i-1}+z_{n-i-1}}{2}=-1-\sqrt{2}+\frac{2\sqrt{2}}{1+\alpha^{n-i}},\\
K&=\frac{-x_{n-i-1}+y_{n-i-1}+z_{n-i-1}}{2}=-1-\sqrt{2}+\frac{2\sqrt{2}}{1-\alpha^{n-i}}.
\end{align*}
\begin{figure}[h!]
\begin{center}
  \psfrag{r}{$p_n$}\psfrag{2}{$v_i$}\psfrag{3}{$q_{i+1}$}\psfrag{4}{$M$} \psfrag{A}{$L_n$}\psfrag{B}{$\Rightarrow$}
   \psfrag{c}{$p_{i}$}\psfrag{1}{$q_{i}$}\psfrag{d}{$u_{n-1}$}\psfrag{2}{$v_i$} \psfrag{e}{$p_{i+1}$}\psfrag{c}{$p_{i}$}
   \psfrag{d}{$u_i$}\psfrag{s}{$N$}\psfrag{t}{$K$} \psfrag{H}{$L_n$} \psfrag{5}{$z_i$}
  \includegraphics[width=60mm]{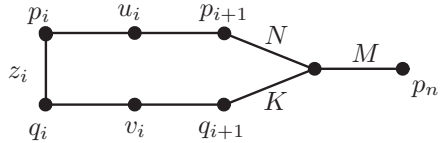}\\
  \caption{$Y$-shaped form of $L_n$.}
\end{center}
\end{figure}

Using parallel and series circuit reductions, we have
\begin{align}
r(p_n,p_i)=&\frac{1}{\frac{1}{N+2}+\frac{1}{K+2+z_i}}+M =\frac{(N+2)(K+2+z_i)}{N+K+z_i+4}+M\notag\\
=&n-i+\frac{(1-\alpha^{n-i})(2-2\alpha^{n+i+2}-\alpha^{n+i+1}-\alpha^{n-i+1}+\alpha^{2i+1}+\alpha)}{4\sqrt{2} (1-\alpha^{2n+2})}.\label{eq:3.6}
\end{align}

Similarly, we can obtain formula for $r(p_n,q_i)$ as
\begin{align}
r(p_n,q_i)=n-i+\frac{(1+\alpha^{n-i})(2+2\alpha^{n+i+2}+\alpha^{n+i+1}+\alpha^{n-i+1}+\alpha^{2i+1}+\alpha)}{4\sqrt{2} (1-\alpha^{2n+2})}.\label{eq:3.7}
\end{align}

On the other hand, $r(p_i,q_i)$ can also be solved (see Fig. 5) as
\begin{align}
r(p_i,q_i)=\frac{1}{\frac{1}{z_i}+\frac{1}{N+K+4}}=\frac{z_i(N+K+4)}{N+K+z_i+4}=\frac{(1+\alpha^{2n-2i+1})(1+\alpha^{2i+1})}{\sqrt{2}(1-\alpha^{2n+2})}.\label{eq:3.8}
\end{align}

Although the formulae in (\ref{eq:3.6})-(\ref{eq:3.8}) are under the condition $n>i>0$, whenever $i=n$ or $i=0$ these formulae are consistent with the result of $x_n$, $y_n$, $z_n$. Therefore, formulae in (\ref{eq:3.6})-(\ref{eq:3.8}) are valid for $n\geq i\geq 0$.

Now we try to obtain formulae for $r(p_i,p_j)$ and $r(q_i,p_j)$ with each integer $n$ and $i$ satisfying $n>i\geq 0$. We can consider $L_n$ as the union of two graphs: One is $L_{n-i-1}$ with two pendant paths and the other is $L_i$. We transform $L_i$ to a $Y$-shaped graph by using Lemma 2.1 and $D$, $E$, $F$ are the resistances along edges in the graph. These reductions are illustrated in Fig. 6, and we obtain the reduced graph in the last stage. By the definition of effective resistance, $D+E=r_{L_i}(p_i,p_j),\, D+F=r_{L_i}(p_i,q_i)=z_i,\, E+F=r_{L_i}(q_i,p_j)$. This gives
\begin{align*}
D=\frac{r_{L_i}(p_i,p_j)+z_i-r_{L_i}(q_i,p_j)}{2},\  E=\frac{r_{L_i}(p_i,p_j)-z_i+r_{L_i}(q_i,p_j)}{2}, \  F= \frac{-r_{L_i}(p_i,p_j)+z_i+r_{L_i}(q_i,p_j)}{2}.
\end{align*}
\begin{figure}[H]
\centering
\psfrag{1}{$p_0$}
\psfrag{2}{$p_j$}
\psfrag{3}{$p_i$}
\psfrag{4}{$u_i$}
\psfrag{5}{$p_{i+1}$}
\psfrag{6}{$p_n$}
\psfrag{a}{$q_0$}
\psfrag{b}{$q_j$}
\psfrag{c}{$q_i$}
\psfrag{d}{$v_i$}
\psfrag{e}{$q_{i+1}$}
\psfrag{f}{$q_n$}
\psfrag{w}{$\dots$}
\psfrag{7}{$\Downarrow$}
\psfrag{8}{$\Rightarrow$}
\psfrag{z}{$z_{n-i-1}$}
\psfrag{x}{$z_{n-i-1}+4$}
\psfrag{r}{$D$}
\psfrag{q}{$E$}
\psfrag{p}{$F$}
\includegraphics[width=4in]{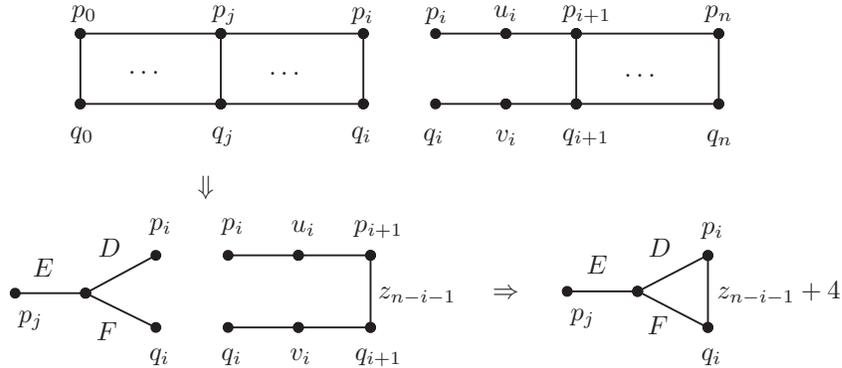}
\caption{Three steps to simplify $L_n$.}
\end{figure}

Using parallel and series circuit reductions, we obtain
\begin{align}
r(p_i,p_j)=&\frac{1}{\frac{1}{D}+\frac{1}{z_{n-i-1}+F+4}}+E=\frac{D(z_{n-i-1}+F+4)}{z_i+z_{n-i-1}+4}+E\notag\\
=&i-j+\frac{1-\alpha^{i-j}}{4\sqrt{2}(1-\alpha^{2n+2})}(2-\alpha^{i+j+1}+\alpha^{2j+1}+\alpha^{2n-2i+1}(1-\alpha^{i-j}-2\alpha^{i+j+1})).\label{eq:3.9}
\end{align}

Similarly, we can obtain formula for $r(q_i,p_j)$ as
\begin{align}
r(q_i,p_j)
=i-j+\frac{1+\alpha^{i-j}}{4\sqrt{2}(1-\alpha^{2n+2})}(2+\alpha^{i+j+1}+\alpha^{2j+1}+\alpha^{2n-2i+1}(1+\alpha^{i-j}+2\alpha^{i+j+1})).\label{eq:3.10}
\end{align}

Let $n=i$ in (\ref{eq:3.9}) and (\ref{eq:3.10}), we can get the formulae which are consistent with (\ref{eq:3.6}) and (\ref{eq:3.7}). So these formulae in (\ref{eq:3.9}) and (\ref{eq:3.10}) are valid for each integers $n$ and $i$ satisfying $0\leq i\leq n$. In other words, we can express the effective resistances between any pair of $p_i$ and $q_j$ in $L_n$.

Our next aim is to obtain the effective resistances from $u_i$ (resp. $v_i$) to each of the rest vertices in $L_n$. First we consider $r(p_n,u_i)$, $r(p_n,v_i)$ and $r(u_i,v_i)$.
\begin{figure}[!ht]
\centering
\psfrag{1}{$p_0$}
\psfrag{2}{$\cdots$}
\psfrag{k}{$v_k$}
\psfrag{b}{$q_0$}
\psfrag{t}{$u_k$}
\psfrag{h}{$L_7'$}
\psfrag{7}{$q_k$}
\psfrag{8}{$R_2$}
\psfrag{u}{$1$}
\psfrag{x}{$p_{k-1}$}
\psfrag{y}{$u_{k-1}$}
\psfrag{z}{$p_k$}
\psfrag{5}{$q_{k-1}$}
\psfrag{6}{$v_{k-1}$}
\psfrag{7}{$q_k$}
  \includegraphics[width=80mm]{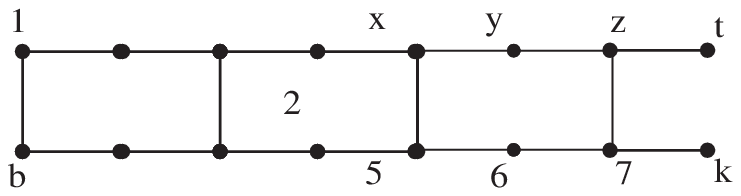}
  \caption{$L_{k}'$.}
\end{figure}

Let $L_{k}'=L_k+\{p_ku_k, q_kv_k\}$ as depicted in Fig. 7. According to the definition of effective resistance, we have
\begin{align*}
a_{k}:=&r_{L_{k}'}(u_{k},p_0)=1+x_{k}=-\sqrt{2}+\frac{2\sqrt{2}}{1+\alpha^{k+1}}+k,\ \ \
b_{k}:=r_{L_{k}'}(u_{k},q_0)=1+y_{k}=-\sqrt{2}+\frac{2\sqrt{2}}{1-\alpha^{k+1}}+k,\\
c_{k}:=&r_{L_{k}'}(u_{k},v_{k})=2+z_{k}=-2\sqrt{2}+\frac{4\sqrt{2}}{1-\alpha^{2k+2}}.
\end{align*}

Then we may consider $L_n$ as the union of two graphs: One is $L_{i+1}'$ and the other is $L_{n-i-2}'$. We transform $L_{n-i-2}'$ to a $Y$-shaped graph by Lemma {\ref{lem1}}, which is depicted in Fig. 8. $X, Y, Z$ are the resistances between vertices of the $Y$-shaped graph. Therefore, we have $X+Z=a_{n-i-2}, Y+Z=b_{n-i-2}, X+Y=c_{n-i-2}$, i.e.,
$$X=-\sqrt{2}+\frac{2\sqrt{2}}{1+\alpha^{n-i-1}},\quad \quad Y=-\sqrt{2}+\frac{2\sqrt{2}}{1-\alpha^{n-i-1}},\quad \quad Z=n-i-2.$$
\begin{figure}[H]
\centering
\psfrag{a}{$u_i$}
\psfrag{b}{$p_{i+1}$}
\psfrag{c}{$u_{i+1}$}
\psfrag{1}{$v_i$}
\psfrag{2}{$q_{i+1}$}
\psfrag{3}{$v_{i+1}$}
\psfrag{d}{$c_i$}
\psfrag{X}{$X$}
\psfrag{Y}{$Y$}
\psfrag{Z}{$Z$}
\psfrag{r}{$\Rightarrow$}
\psfrag{p}{$p_n$}
\psfrag{4}{$X+1$}
\psfrag{5}{$Y+1$}
\psfrag{6}{$Z$}
\psfrag{7}{$1$}
\psfrag{8}{$c_i+1$}
\psfrag{R}{$X'$}
\psfrag{S}{$Y'$}
\psfrag{T}{$Z'$}
\includegraphics[width=6in]{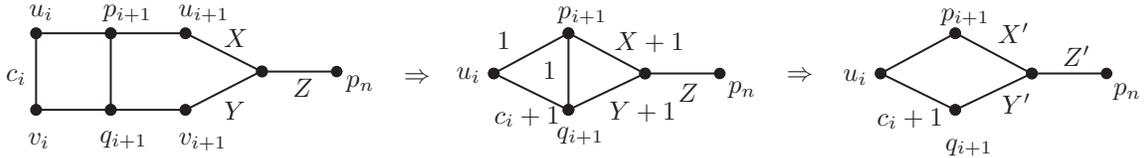}
\caption{Three steps to simplify $L_n$.}
\end{figure}

We use parallel and series circuit reductions to obtain
$$r(u_i,v_i)=\frac{1}{\frac{1}{c_i}+\frac{1}{2+\frac{1}{1+\frac{1}{X+Y+2}}}}=\frac{\sqrt{2}(1+\alpha^{2i+2})(1+\alpha^{2n-2i})}{1-\alpha^{2n+2}}.$$
Using Lemma {\ref{lem1}} and transform the $\Delta$-shaped graph to the $Y$-shaped graph, i.e., the last graph in Fig. 8, where $X'=\frac{X+1}{X+Y+3}$, $Y'=\frac{Y+1}{X+Y+3}$ and $Z'=Z+\frac{(X+1)(Y+1)}{X+Y+3}$. Then we may figure out the formulae for $r(p_n,u_i)$, $r(p_n,v_i)$ and $r(u_i,v_i)$ as
\begin{align*}
r(p_n,u_i)=&\frac{1}{\frac{1}{X'+1}+\frac{1}{Y'+c_i+1}}+Z'
=\frac{\sqrt{2}(1-\alpha^{n-i})(3+\alpha^{2i+2}-\alpha^{n-i}-3\alpha^{n+i+2})}{4(1-\alpha^{2n+2})}+n-i-1,\\[5pt]
r(p_n,v_i)=&\frac{\sqrt{2}(1+\alpha^{n-i})(3+\alpha^{2i+2}+\alpha^{n-i}+3\alpha^{n+i+2})}{4(1-\alpha^{2n+2})}+n-i-1.
\end{align*}
\begin{figure}[!ht]
\centering
\psfrag{r}{$p_n$}
\psfrag{a}{$p_j$}
\psfrag{4}{$T$}
\psfrag{s}{$R$}
\psfrag{t}{$S$}
\psfrag{2}{$\cdots$}
\psfrag{k}{$n-1$}
\psfrag{f}{$L_n$}
\psfrag{h}{$\Rightarrow$}
\psfrag{l}{$c_{n-i-1}$}
\psfrag{m}{$q_n$}
\psfrag{r}{$p_n$}
\psfrag{x}{$p_i$}
\psfrag{y}{$u_i$}
\psfrag{z}{$p_{i+1}$}
\psfrag{5}{$q_i$}
\psfrag{6}{$v_i$}
\psfrag{7}{$q_{i+1}$}
  \includegraphics[width=120mm]{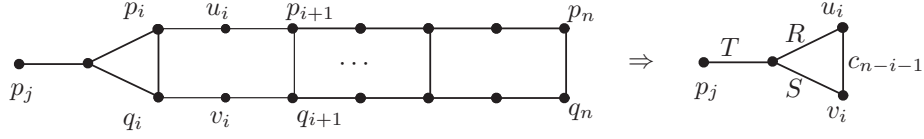}
  \caption{$Y$-shaped $L_n$ with effective resistance $R, S, T$.}
\end{figure}

Next, we are to obtain formulae for $r(u_i,p_j)$ and $r(v_i,p_j)$. We also use Lemma {\ref{lem1}} and consider $L_n$ as the union of $L_i'$ and $L_{n-i-1}'$. We transform $L_i'$ to a $Y$-shaped graph and $L_{n-i-1}'$ to an edge with value $c_{n-i-1}$. These reductions are illustrated in Fig. 9. By the definition of effective resistance, $R+T=1+r_{L_i}(p_i,p_j)$, $S+T=1+r_{L_i}(p_i,q_j)$ and $R+S=2+r_{L_i}(p_i,q_i)=2+z_i$. For each $0 \leq j \leq i\leq n$, we have
\begin{align}\label{eq:3.11}
r(u_i,p_j)=&\frac{1}{\frac{1}{c_{n-i-1}+S}+\frac{1}{R}}+T=i-j+\frac{f(\alpha)}{4\sqrt{2}(1-\alpha^{2n+2})},
\end{align}
where $f(\alpha)=2\alpha^{2i+2}+\alpha^{2n-2j+1}+2\alpha^{2n-2i}+\alpha^{2n+3}-\alpha^{-i-j}(\alpha+1)(1+\alpha^{2j+1})(\alpha^{2i}+\alpha^{2n})+\alpha^{-1}+\alpha^{2j+1}.$

Similarly,
\begin{align}\label{eq:3.12}
r(v_i,p_j)=&\frac{1}{\frac{1}{c_{n-i-1}+R}+\frac{1}{S}}+T=i-j+\frac{g(\alpha)}{4\sqrt{2}(1-\alpha^{2n+2})},
\end{align}
where $g(\alpha)=2\alpha^{2i+2}+\alpha^{2n-2j+1}+2\alpha^{2n-2i}+\alpha^{2n+3}+\alpha^{-i-j}(\alpha+1)(1+\alpha^{2j+1})(\alpha^{2i}+\alpha^{2n})+\alpha^{-1}+\alpha^{2j+1}.$
\begin{figure}[!ht]
\centering
\psfrag{r}{$p_n$}
\psfrag{a}{$u_j$}
\psfrag{4}{$W$}
\psfrag{s}{$V$}
\psfrag{t}{$U$}
\psfrag{2}{$\cdots$}
\psfrag{k}{$n-1$}
\psfrag{f}{$L_n$}
\psfrag{h}{$\Rightarrow$}
\psfrag{l}{$c_{n-i-1}$}
\psfrag{m}{$q_n$}
\psfrag{r}{$p_n$}
\psfrag{x}{$p_i$}
\psfrag{y}{$u_i$}
\psfrag{z}{$p_{i+1}$}
\psfrag{5}{$q_i$}
\psfrag{6}{$v_i$}
\psfrag{7}{$q_{i+1}$}
  \includegraphics[width=120mm]{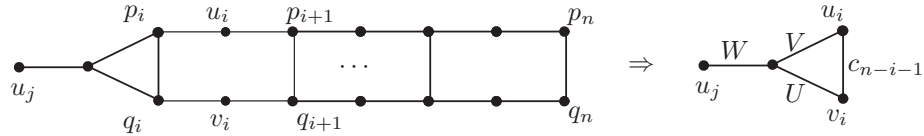}
  \caption{$Y$-shaped $L_n$ with effective resistances $U, V$ and $W$.}
\end{figure}

Finally, using the above results, we obtain the formulae for $r(u_i,u_j)$ and $r(v_i,u_j)$. These reductions are illustrated in Fig. 10, where $U, V, W$ satisfy $U+W=1+r_{L_i}(p_i,u_j)$, $V+W=1+r_{L_i}(p_i,v_j)$ and $U+V=2+r_{L_i}(p_i,q_i)=2+z_i$. For each $0 \leq j \leq i\leq n$, we have
\begin{align}
r(u_i,u_j)=&\frac{1}{\frac{1}{c_{n-i-1}+V}+\frac{1}{U}}+W\notag\\
=&i-j+\frac{(\alpha^{i+1}-\alpha^{j+1})(\alpha^{i+1}-\alpha^{j+1}-2\alpha^{-j-1}+\alpha^{2n-i-2j-1}-\alpha^{2n-j-2i-1}+2\alpha^{2n-i+1})}{2\sqrt{2}(1-\alpha^{2n+2})}.\label{eq:3.13}
\end{align}
Similarly,
\begin{align}\label{eq:3.14}
r(v_i,u_j)=i-j+\frac{(\alpha^{i+1}+\alpha^{j+1})(\alpha^{i+1}+\alpha^{j+1}+2\alpha^{-j-1}+\alpha^{2n-i-2j-1}+\alpha^{2n-j-2i-1}+2\alpha^{2n-i+1})}{2\sqrt{2}(1-\alpha^{2n+2})}.
\end{align}
\subsection{\normalsize The maximum and the minimum effective resistances in $L_n$}

In this subsection, we consider the extremal problems on the effective resistances in $L_n$. 
\begin{lem}\label{lem4.1}
Assume $0 \leq j \leq n-k$ with $k>0$. Let $L_n$ be the linear hexagonal chain as depicted in Fig. $1$.
\begin{wst}
\item[{\rm (i)}] For any fixed $k, \, r(p_j,p_{j+k})$ is convex in $i$, i.e.,
{\begin{align*}
r&(p_0,p_{k})>r(p_1,p_{1+k})>\dots >r(p_{\lfloor \frac{n-k}{2} \rfloor},p_{\lfloor \frac{n-k}{2} \rfloor+k})
=r(p_{\lceil \frac{n-k}{2} \rceil},p_{\lceil \frac{n-k}{2} \rceil+k}) <\dots<r(p_{n-k},p_n).
\end{align*}}
\item[{\rm (ii)}] For any fixed j, $r(p_j,p_{j+k})$ is monotone increasing in $k$.
\item[{\rm (iii)}] {$r(p_{\lfloor \frac{n-k}{2} \rfloor},p_{\lfloor \frac{n-k}{2} \rfloor+k})$} is monotone increasing in $k$.
\end{wst}
\end{lem}
\begin{proof}
(i)\ For convenience, let $i:=j+k$ and consider $r(p_{j+1},p_{i+1})-r(p_j,p_i)$.

According to the formulae ({\ref{eq:3.9}}) in section 3, we obtain that
$${r(p_{j+1},p_{i+1})-r(p_j,p_i)=\frac{\alpha^{-1-2i-2j}(1-\alpha^{2})(\alpha^{i}-\alpha^{j})^{2}(\alpha^{2n}-\alpha^{2i+2j+2})}{4\sqrt{2} (1-\alpha^{2n+2})}.}$$
Hence, $r(p_{j+1},p_{i+1})  \geq  r(p_j,p_i)$ if $n\leq i+j+2$ and $r(p_{j+1},p_{i+1})  <  r(p_j,p_i)$ otherwise. So (i) is proved.

(ii)\ It suffices to show that $r(p_j,p_{i+1})-r(p_j,p_i)>0 $ for fixed $i$ and $j$. In fact,
{\begin{align*}
r(p_j,p_{i+1})-r(p_j,p_i)&=1+\frac{1}{4\sqrt{2}(1-\alpha^{2n+2})}(1-\alpha)\alpha^{-1-2i-j}f(\alpha,i,j),
\end{align*}}
where $f(\alpha,i,j)= 2\alpha^{3i+1}-\alpha^{4i+j+2}-\alpha^{4i+j+3}+2\alpha^{3i+2j+2}-2\alpha^{2n+2j+i+2}+\alpha^{2n+j}+\alpha^{2n+j+1}-2\alpha^{2n+i+1}.$

Note that $2\alpha^{3i+1}-\alpha^{4i+j+2}-\alpha^{4i+j+3}>0, 2\alpha^{3i+2j+2}-2\alpha^{2n+2j+i+2}>0$ and $\alpha^{2n+j}+\alpha^{2n+j+1}-2\alpha^{2n+i+1}>0$. Hence, we have
$r(p_{j+1},p_i)>r(p_j,p_i),$ as desired.

(iii) Based on (\ref{eq:3.9}), it suffices to consider the difference
\[\label{eq:4.1}
  r(p_{\lfloor \frac{n-k-1}{2} \rfloor},p_{\lfloor \frac{n-k-1}{2} \rfloor+k+1})-r(p_{\lfloor \frac{n-k}{2} \rfloor},p_{\lfloor \frac{n-k}{2} \rfloor+k}).
\]
If $n-k$ is even, then the above difference equals to
$$
 1+\frac{1-\alpha}{4\sqrt{2}(1-\alpha^{2n+2})}(\alpha^{-1-k+n}(1+\alpha)(1-\alpha^{2+2k})+2(\alpha^{k}-\alpha^{n}+\alpha^{1+n}-\alpha^{1-k+2n})).
$$
Note that {$1-\alpha^{2+2k}> 0$, $\alpha^{k}-\alpha^{n}+\alpha^{1+n}-\alpha^{1-k+2n}>0$ and $\frac{1-\alpha}{4\sqrt{2}(1-\alpha^{2n+2})}>0$}. Hence, the difference in (\ref{eq:4.1}) is positive.

If $n-k-1$ is odd, then the above difference (\ref{eq:4.1}) equals to
$$
1+\frac{1-\alpha}{4\sqrt{2}(1-\alpha^{2n+2})}(\alpha^{-k+n}(1+\alpha)(1-\alpha^{2k})+2(\alpha^{k}+\alpha^{n}-\alpha^{1+n}-\alpha^{1+2n})).
$$
As {$1-\alpha^{2k}>0$, $\alpha^{k}+\alpha^{n}-\alpha^{1+n}-\alpha^{1+2n}>0
$ and $\frac{1-\alpha}{4\sqrt{2}(1-\alpha^{2n+2})}>0$}, we have the difference in (\ref{eq:4.1}) is positive. Hence, (iii) holds.
\end{proof}

Similarly, we can obtain this property for $r(q_j,p_i)$, $r(u_j,u_i)$ and $r(v_j,u_i)$ as
\begin{lem}\label{lem4.2}
Assume $0 \leq j \leq n-k$ with $k\geq 0$. Let $L_n$ be the linear hexagonal chain as depicted in Fig. $1$.
\begin{wst}
\item[{\rm (i)}] For any fixed $k$, $\{r(q_j,p_{j+k})\}_{j=0}^{n-k}\, ($resp. $\{r(u_j,u_{j+k})\}_{j=0}^{n-k},\, \{r(v_j,u_{j+k})\}_{j=0}^{n-k})$ is convex in $j$, i.e.,
{\begin{align*}
r&(q_0,p_{k})>r(q_1,p_{1+k})>\dots >r(q_{\lfloor \frac{n-k}{2} \rfloor},p_{\lfloor \frac{n-k}{2} \rfloor+k})
=r(q_{\lceil \frac{n-k}{2} \rceil},p_{\lceil \frac{n-k}{2} \rceil+k})<\dots<r(q_{n-k},p_n).\\
r&(u_0,u_{k})>r(u_1,u_{1+k})>\dots >r(u_{\lfloor \frac{n-k-1}{2} \rfloor},u_{\lfloor \frac{n-k-1}{2} \rfloor+k})
=r(u_{\lceil \frac{n-k}{2} \rceil},u_{\lceil \frac{n-k}{2} \rceil+k})<\dots<r(u_{n-k},u_n),\, k> 0.\\
r&(v_0,u_{k})>r(v_1,u_{1+k})>\dots >r(v_{\lfloor \frac{n-k}{2} \rfloor},u_{\lfloor \frac{n-k}{2} \rfloor+k})
=r(v_{\lceil \frac{n-k}{2} \rceil},u_{\lceil \frac{n-k}{2} \rceil+k}) <\dots<r(v_{n-k},u_n).
\end{align*}}
\item[{\rm (ii)}]For any fixed $j$,  $r(q_j,p_{j+k}), r(u_j,u_{j+k})$ and $r(v_j,u_{j+k})$ are monotone increasing in $k$, respectively.
\item[{\rm (iii)}] {$r(q_{\lfloor \frac{n-k}{2} \rfloor},p_{\lfloor \frac{n-k}{2} \rfloor+k})$, $r(u_{\lfloor \frac{n-k}{2} \rfloor},u_{\lfloor \frac{n-k}{2} \rfloor+k})$ and $r(v_{\lfloor \frac{n-k}{2} \rfloor},u_{\lfloor \frac{n-k}{2} \rfloor+k})$} are monotone increasing in $k$, respectively.
\end{wst}
\end{lem}
\begin{lem}\label{lem4.3}
Assume $0 \leq j \leq n-k$ with $k\geq 0$. Let $L_n$ be the linear hexagonal chain as depicted in Fig. $1$.
\begin{wst}
\item[{\rm (i)}]For any fixed $k, \, r(p_j,u_{j+k})$ and $r(p_j,v_{j+k})$ are convex in $i$, i.e.,
{\begin{align*}
r&(p_0,u_{k})>r(p_1,u_{1+k})>\dots >r(p_{\lfloor \frac{n-k}{2} \rfloor},u_{\lfloor \frac{n-k}{2} \rfloor+k})
=r(p_{\lceil \frac{n-k}{2} \rceil},u_{\lceil \frac{n-k}{2} \rceil+k})
<\dots<r(p_{n-k},u_n),\\
r&(p_0,v_{k})>r(p_1,v_{1+k})>\dots >r(p_{\lfloor \frac{n-k}{2} \rfloor},v_{\lfloor \frac{n-k}{2} \rfloor+k})
=r(p_{\lceil \frac{n-k}{2} \rceil},v_{\lceil \frac{n-k}{2} \rceil+k})<\dots<r(p_{n-k},v_n).
\end{align*}}
\item[{\rm (ii)}] For any fixed $j,\, r(p_j,u_{j+k})$ and $r(p_j,v_{j+k})$ are monotone increasing in $k$.
\item[{\rm (iii)}] {$r(p_{\lfloor \frac{n-k-1}{2} \rfloor},u_{\lfloor \frac{n-k-1}{2} \rfloor+k})$ and $r(p_{\lfloor \frac{n-k-1}{2} \rfloor},v_{\lfloor \frac{n-k-1}{2} \rfloor+k})$} are monotone increasing in $k$.
\end{wst}
\end{lem}
\begin{proof}
By the symmetry, we only show the proof for $r(p_j,u_{j+k})$ in what follows. We omit the procedure for $r(p_j,v_{j+k})$.

(i) Let $i:=j+k$. Then 
\begin{align*}
r(p_{j+1},u_{i+1})-r(p_j,u_i)=&\frac{1}{4\sqrt{2}(1-\alpha^{2n+2})}\alpha^{-2(i+j+1)}(1-\alpha^{2})\left[\alpha^{2n}(\alpha^{1+2i}+2\alpha^{2j}-\alpha^{i+j}(1+\alpha))\right.\\
                              &-\left.\alpha^{2i+2j+3}(2\alpha^{1+2i}+\alpha^{2j}-\alpha^{i+j}(1+\alpha))\right].
\end{align*}
It's easy to see that $\frac{1}{4\sqrt{2}(1-\alpha^{2n+2})}\alpha^{-2(i+j+1)}(1-\alpha^{2})>0$.

If $n\geq i+j+2$, we have
\begin{align*}
&\alpha^{2n}(\alpha^{1+2i}+2\alpha^{2j}-\alpha^{i+j}(1+\alpha))-\alpha^{2i+2j+3}(2\alpha^{1+2i}+\alpha^{2j}-\alpha^{i+j}(1+\alpha))\\
&<(\alpha^{2n}-\alpha^{2i+2j+3})(\alpha^{1+2i}+2\alpha^{2j}-\alpha^{i+j}(1+\alpha))\tag{As $\alpha^{1+2i}<\alpha^{2j}$}\\
&<0.
\end{align*}

If $n<i+j+2$,
\begin{align*}
&\alpha^{2n}(\alpha^{1+2i}+2\alpha^{2j}-\alpha^{i+j}(1+\alpha))-\alpha^{2i+2j+3}(2\alpha^{1+2i}+\alpha^{2j}-\alpha^{i+j}(1+\alpha))\\
&=(\alpha^{2n}-\alpha^{2i+2j+3})(\alpha^{1+2i}+2\alpha^{2j}-\alpha^{i+j}(1+\alpha))+\alpha^{2i+4j+3}-\alpha^{4i+2j+4}\\
&>0.
\end{align*}
The last inequality is due to $\alpha^{2n}-\alpha^{2i+2j+3}\geq 0$, $\alpha^{1+2i}+2\alpha^{2j}-\alpha^{i+j}(1+\alpha)>0$ and $\alpha^{2i+4j+3}-\alpha^{4i+2j+4}>0$.
So (i) is proved.

(ii) We show that $r(p_j,u_{i+1})-r(p_j,u_i)>0 $ for fixed $i$ and $j$.
\begin{align*}
r(p_j,u_{i+1})-r(p_j,u_i)=&\ 1+\frac{1-\alpha^{2}}{4\sqrt{2}(1-\alpha^{2n+2})}\left[\alpha^{i-j}(1-2\alpha^{i+j+2}+\alpha^{2j+1})\right.\\
                           &\ \left.+\alpha^{2n-i-j-2}(-\alpha+2\alpha^{j-i}-\alpha^{2j+2})\right].
\end{align*}
Note that $\frac{1-\alpha^{2}}{4\sqrt{2}(1-\alpha^{2n+2})}>0$, $1-2\alpha^{i+j+2}+\alpha^{2j+1}>0$ and $-\alpha+2\alpha^{j-i}-\alpha^{2j+2}>0$. Hence, we have
$r(p_j,u_{i+1})>r(p_j,u_i),$ as desired.

(iii) It follows directly by (\ref{eq:3.11}) and (\ref{eq:3.12}). 
\end{proof}


The next lemma follows directly from Lemmas \ref{lem4.1}-\ref{lem4.3}.
\begin{thm}\label{thm4.4}
For the graph $L_{n}$ with $n\geq1$, we have
\begin{wst}
\item[{\rm (i)}] $r(p_0,p_{n})\geq r(p_j,p_i)\geq r(p_{\lfloor \frac{n-1}{2} \rfloor},p_{\lfloor \frac{n-1}{2} \rfloor+1})$ for $0\leq j < i \leq n$.
\item[{\rm (ii)}]$r(q_0,p_{n})\geq r(q_j,p_i)\geq r(q_{\lfloor \frac{n}{2} \rfloor},p_{\lfloor \frac{n}{2} \rfloor})$ for $0\leq j \leq i \leq n$.
\item[{\rm (iii)}] $r(u_0,u_{n-1})\geq r(u_j,u_i)\geq r(u_{\lfloor \frac{n-1}{2} \rfloor},u_{\lfloor \frac{n-1}{2} \rfloor+1})$ for $0\leq j < i \leq n-1$.
\item[{\rm (iv)}] $r(v_0,u_{n-1})\geq r(v_j,u_i)\geq r(v_{\lfloor \frac{n}{2} \rfloor},u_{\lfloor \frac{n}{2} \rfloor})$ for $0\leq j \leq i \leq n-1$.
\item[{\rm (v)}] $r(p_0,u_{n-1})\geq r(p_j,u_i)\geq r(p_{\lfloor \frac{n}{2} \rfloor},u_{\lfloor \frac{n}{2} \rfloor})$ for $0\leq j \leq i \leq n-1$.
\item[{\rm (vi)}] $r(p_0,v_{n-1})\geq r(p_j,v_i)\geq r(p_{\lfloor \frac{n}{2} \rfloor},v_{\lfloor \frac{n}{2} \rfloor})$ for $0\leq j \leq i \leq n-1$.
\end{wst}
\end{thm}

Based on Theorem \ref{thm4.4}, we may determine the maximum and the minimum effective resistances in $L_{n}$.
\begin{thm}
For the graph $L_{n}$ with $n\geq1$ and any vertex $a,b\in V(L_{n})$, we have
$r(p_0,q_{n})\geq r(a,b)\geq r(p_{\lfloor \frac{n}{2} \rfloor},q_{\lfloor \frac{n}{2} \rfloor}).$
\end{thm}
\begin{proof}
First, we proof $r(p_0,q_n)\geq r(a,b)$.
Note that
$$
\begin{array}{ll}
 r(p_0,p_{n})=n-1-\sqrt{2}+\frac{2\sqrt{2}}{1+\alpha^{n+1}},  & r(p_0,q_{n})=n-1-\sqrt{2}+\frac{2\sqrt{2}}{1-\alpha^{n+1}}, \\
  r(p_0,u_{n-1})=n-7-5\sqrt{2}+\frac{2\sqrt{2}}{1-\alpha^{n+1}}+\frac{8\sqrt{2}}{1+\alpha^{n+1}}, & r(p_0,v_{n-1})=n-7-5\sqrt{2}+\frac{8\sqrt{2}}{1-\alpha^{n+1}}+\frac{2\sqrt{2}}{1+\alpha^{n+1}}, \\
  r(u_0,u_{n-1})=n-13-9\sqrt{2}+\frac{18\sqrt{2}}{1+\alpha^{n+1}}, &  r(u_0,v_{n-1})=n-13-9\sqrt{2}+\frac{18\sqrt{2}}{1-\alpha^{n+1}}.
\end{array}
$$
This gives, for $n\geq 1$, that
\begin{align*}
r(p_0,q_{n})-r(p_0,p_{n})=&\frac{4\alpha^{n}(-4+3\sqrt{2})}{1-\alpha^{2n+2}}> 0,\\
r(p_0,q_{n})-r(p_0,u_{n-1})=&6+4\sqrt{2}-\frac{8\sqrt{2}}{1+\alpha^{n+1}}>0,\\
r(p_0,q_{n})-r(p_0,v_{n-1})=&6+4\sqrt{2}-\frac{6\sqrt{2}}{1-\alpha^{n+1}}-\frac{2\sqrt{2}}{1+\alpha^{n+1}}>0,\\
r(p_0,q_{n})-r(u_0,u_{n-1})=&12+8\sqrt{2}+\frac{2\sqrt{2}}{1-\alpha^{n+1}}-\frac{18\sqrt{2}}{1+\alpha^{n+1}}>0,\\
r(p_0,q_{n})-r(u_0,v_{n-1})=&12+8\sqrt{2}-\frac{16\sqrt{2}}{1-\alpha^{n+1}}>0.
\end{align*}
By Theorem \ref{thm4.4}, we get $r(p_0,q_n)\geq r(a,b)$.

Similarly, we may show that $r(a,b)\geq r(p_{\lfloor \frac{n}{2} \rfloor},q_{\lfloor \frac{n}{2} \rfloor})$ according to the parity of $ n $. We omit the procedure here.
\end{proof}

At the end of this section, we turn to the asymptotic properties of resistance distances in $L_n$.
\begin{thm}
\begin{wst}
\item[{\rm (i)}] For all fixed $i$ and $j$, one has
\begin{align*}
  \lim_{n\to \infty}r(p_i,p_j)&=i-j+\frac{(1-\alpha^{i-j})(2-\alpha^{i+j+1}+\alpha^{2j+1})}{4\sqrt{2}},\\
  \lim_{n\to \infty}r(q_i,p_j)&=i-j+\frac{(1+\alpha^{i-j})(2+\alpha^{i+j+1}+\alpha^{2j+1})}{4\sqrt{2}}.
\end{align*}
\item[{\rm (ii)}] $\lim\limits_{n \to \infty}\frac{1}{n}r(p_n,p_0)=\lim\limits_{n \to \infty}\frac{1}{n}r(q_n,p_0)=1.$
\item[{\rm (iii)}]
  $\lim\limits_{n \to \infty} [r_{L_{n+1}}(p_{n+1},p_0)-r_{L_n}(p_n,p_0)]=
  \lim\limits_{n \to \infty} [r_{L_{n+1}}(q_{n+1},p_0)-r_{L_n}(q_n,p_0)]=1$
\item[{\rm (iv)}]
  $\lim\limits_{n \to \infty}r(p_{\lfloor \frac{n-1}{2} \rfloor+1},p_{\lfloor \frac{n-1}{2} \rfloor})=2-\frac{\sqrt{2}}{2},
  \quad \lim\limits_{n \to \infty}r(p_{\lfloor \frac{n}{2} \rfloor},q_{\lfloor \frac{n}{2} \rfloor})=\frac{\sqrt{2}}{2}$
\end{wst}
\end{thm}
\begin{proof}
(i)-(iv) follow directly by ({\ref{eq:3.9}}) and (\ref{eq:3.10}).
\end{proof}
By a similar discussion, it is not difficult to determine the limit value on the resistance distance between any other pair of vertices in $L_n$ as $n\rightarrow \infty$. Here we omit the contents.
\section{\normalsize The effective resistance in $R_n$}\setcounter{equation}{0}
\subsection{\normalsize Determining the effective resistance between any two vertices in $R_n$}

In this subsection, we determine the resistance distance for every pair of vertices of $R_n$. We abbreviate $r_{R_n}(u,v)$ to $r(u,v)$ in this section for any pair of vertices $u,v$ in $R_n$. First of all, we need the following lemma, which simplifies the circuit of $R_n$.
\begin{figure}[H]
\centering
\psfrag{u}{$u$}
\psfrag{w}{$w$}
\psfrag{t}{$t$}
\psfrag{2}{$\cdots$}
\psfrag{h}{$\Rightarrow$}
\psfrag{m}{$q_n$}
\psfrag{r}{$p_n$}
\psfrag{x}{$p_0$}
\psfrag{y}{$u_0$}
\psfrag{z}{$p_{1}$}
\psfrag{5}{$q_0$}
\psfrag{6}{$v_0$}
\psfrag{7}{$q_{1}$}
\includegraphics[width=120mm]{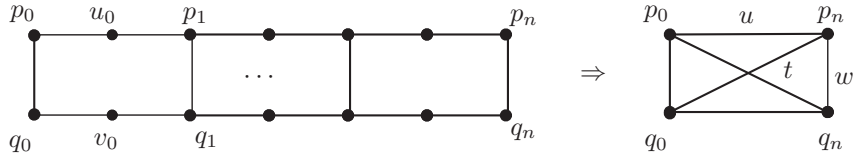}
\caption{Simplified circuit of $L_n$.}
\end{figure}
\begin{lem}\label{lem:5.1}
Assume that $L_n$ is a linear hexagonal chain, then we can transform $L_n$ to $K_4$, where $w$, $u$ and $t$ are the resistances along edges in $K_4$ (see Fig. 11). This transformation yield equivalent circuits.
\end{lem}
\begin{proof}
To prove this lemma, we just need to prove the existence of such resistances $w$, $u$ and $t$ satisfying the following equations:
\begin{align}
r_{K_4}(p_n,q_n)=&r_{K_4}(p_0,q_0)=f(w,u,t)=r_{L_n}(p_n,q_n)=r_{L_n}(p_0,q_0),\label{eq:5.2}\\
r_{K_4}(p_n,p_0)=&r_{K_4}(q_n,q_0)=g(w,u,t)=r_{L_n}(p_n,p_0)=r_{L_n}(q_n,q_0),\label{eq:5.3}\\
r_{K_4}(p_n,q_0)=&r_{K_4}(q_n,p_0)=h(w,u,t)=r_{L_n}(p_n,q_0)=r_{L_n}(q_n,p_0),\label{eq:5.4}
\end{align}
where $f(w,u,t)$, $g(w,u,t)$ and $h(w,u,t)$ are functions of $w$, $u$ and $t$.

Our aim is to obtain the formulae of $f(w,u,t)$, $g(w,u,t)$ and $h(w,u,t)$. With the ordering of the vertices $(p_n,q_n,p_1,q_1)$, the discrete Laplacian matrix $M$ of the graph $K_4$ is as follows:
$$
M=\left(
  \begin{array}{rrrrrrrrr}
    \frac{1}{t}+\frac{1}{u}+\frac{1}{w} & -\frac{1}{w} & -\frac{1}{u} & -\frac{1}{t}\\[2pt]
    -\frac{1}{w} & \frac{1}{t}+\frac{1}{u}+\frac{1}{w} & -\frac{1}{t} & -\frac{1}{u}\\[2pt]
    -\frac{1}{u} & -\frac{1}{t} & \frac{1}{t}+\frac{1}{u}+\frac{1}{w} & -\frac{1}{w}\\[2pt]
    -\frac{1}{t} & -\frac{1}{u} & -\frac{1}{w} & \frac{1}{t}+\frac{1}{u}+\frac{1}{w}
  \end{array}
 \right).
$$
Then we obtain the Moore-Penrose inverse $M^+$ of $M$ (see \cite{Bapat3}) as
$$M^{+}=(L+\frac{1}{4}J)^{-1}+\frac{1}{4}J,$$
where $J$ is the matrix with all entries 1.

We use the Moore-Penrose inverse and Lemma \ref{lem2} to obtain the following
$$f(w,u,t)=\frac{1}{2}(\frac{wt}{w+t}+\frac{uw}{u+w}),\quad \quad g(w,u,t)=\frac{1}{2}(\frac{uw}{u+w}+\frac{ut}{u+t}),\quad \quad h(w,u,t)=\frac{1}{2}(\frac{ut}{u+t}+\frac{wt}{w+t}).$$
Then we have
\begin{align}
\frac{wt}{w+t}&=r_{L_n}(p_n,q_n)+r_{L_n}(p_n,q_{0})-r_{L_n}(p_n,p_{0}):=A_n,\label{eq:5.5}\\
\frac{uw}{u+w}&=r_{L_n}(p_n,q_n)-r_{L_n}(p_n,q_{0})+r_{L_n}(p_n,p_{0}):=B_n,\label{eq:5.6}\\
\frac{ut}{u+t}&=-r_{L_n}(p_n,q_n)+r_{L_n}(p_n,q_{0})+r_{L_n}(p_n,p_{0}):=C_n.\label{eq:5.7}
\end{align}

Based on (3.1)-(3.3), we obtain the formulae of $A_n$, $B_n$ and $C_n$ as
\begin{align*}
A_n=-2-2\sqrt{2}+\frac{4\sqrt{2}}{1-\alpha^{n+1}},\quad \quad B_n=-2-2\sqrt{2}+\frac{4\sqrt{2}}{1+\alpha^{n+1}},
\quad \quad C_n=2n.
\end{align*}
Thus,
\begin{align*}
w=\frac{2}{\frac{1}{A_n}+\frac{1}{B_n}-\frac{1}{C_n}}&=\frac{4(1-\sqrt{2}+(1+\sqrt{2})\alpha^{n+1})(-1+\sqrt{2}+(1+\sqrt{2})\alpha^{n+1})n}{1-2(-1+\sqrt{2})(n+1)-(1+2(1+\sqrt{2})n\alpha^{2n+2})}>0,\\
u=\frac{2}{-\frac{1}{A_n}+\frac{1}{B_n}+\frac{1}{C_n}}&=\frac{4(1-\sqrt{2}+(1+\sqrt{2})\alpha^{n+1})(-1+\sqrt{2}+(1+\sqrt{2})\alpha^{n+1})n}{-3+2\sqrt{2}+(3+2\sqrt{2})\alpha^{2n+2}-4\sqrt{2}n\alpha^{n+1}}>0,\\
t=\frac{2}{-\frac{1}{A_n}+\frac{1}{B_n}+\frac{1}{C_n}}&=\frac{4(1-\sqrt{2}+(1+\sqrt{2})\alpha^{n+1})(-1+\sqrt{2}+(1+\sqrt{2})\alpha^{n+1})n}{-3+2\sqrt{2}+(3+2\sqrt{2})\alpha^{2n+2}+4\sqrt{2}n\alpha^{n+1}}>0.
\end{align*}
Hence, (\ref{eq:5.2}), (\ref{eq:5.3}) and (\ref{eq:5.4}) hold.
\end{proof}

Then we are to determine the resistance distance for every pair of vertices of $R_n$. For convenience, we label the $i$th 6-cycle of $R_n$ as $p_{i-1}, u_{i-1}, p_i, q_i, v_{i-1}, q_{i-1}$ for $i=1,2,\ldots, n$ and let $p_0=p_n, q_0=q_n, u_0=u_n, v_0=v_n$ (see also Fig. 12).  Clearly, the order of $R_n$ is $4n$ and its size is $5n$. According to symmetry of $R_n$, for $i, j \in [n]$, it can be easy to see that
$$
\begin{array}{cc}
  r(p_{k},p_j)=r(q_{k},q_j)=r(p_1,p_i),  & r(p_{k},q_j)=r(q_{k},p_j)=r(p_1,q_i), \\[5pt]
   r(p_{k},u_j)=r(q_{k},v_j)=r(p_1,u_i),  & r(p_{k},v_j)=r(q_{k},u_j)=r(p_1,v_i), \\[5pt]
   r(u_{k},u_j)=r(v_{k},v_j)=r(u_1,u_i),   & r(u_{k},v_j)=r(v_{k},u_j)=r(u_1,v_i),
\end{array}
$$
where $1\leq k\leq n$ and $k\equiv j-i\pmod{n}$.

It suffices for us to determine $r(p_1,p_i), r(p_1,q_i), r(p_1,u_i), r(p_1,v_i), r(u_1,u_i)$ and $r(u_1,v_i)$, respectively.
We firstly determine the formulae for $r(p_1,p_i)$ and $r(p_1,q_i)$ by Moore-Penrose inverse.

Now we use Lemma \ref{lem:5.1} to simplify the circuit of $R_n$. We transform the graph $R_n$ to $R_n'$. Note that we only consider the vertices $p_1$, $p_i$ and $q_i$ in this case, so we replace the path $p_1u_np_n$ (resp. $q_1v_nq_n$, $p_{i-1}u_{i-1}p_{i}$ and $q_{i-1}v_{i-1}q_{i}$) with the edge $p_1p_n$ (resp. $q_1q_n$, $p_{i-1}p_{i}$ and $q_{i-1}q_{i}$) with resistance of two. Then we obtain the graph $R_n''$ (see also Fig. 12).
\begin{figure}[!ht]
\centering
\psfrag{r}{$i$}\psfrag{a}{$a$}\psfrag{b}{$b$}\psfrag{c}{$c$}
\psfrag{2}{$L_{n-i}$}\psfrag{s}{$s$}\psfrag{m}{$m$}\psfrag{w}{$k$}
\psfrag{k}{$p_n$}\psfrag{e}{$R_n''$}
\psfrag{f}{$R_n$}
\psfrag{g}{$R_n'$}
\psfrag{t}{$p_1$}
\psfrag{h}{$q_1$}
\psfrag{7}{$q_i$}
\psfrag{8}{$2$}
\psfrag{u}{$q_n$}
\psfrag{x}{$p_{i-1}$}
\psfrag{y}{$u_{i-1}$}
\psfrag{z}{$p_i$}
\psfrag{5}{$q_{i-1}$}
\psfrag{6}{$v_{i-1}$}
\psfrag{1}{$L_{i-2}$}
\psfrag{0}{$u_n$}
\psfrag{9}{$v_n$}
  \includegraphics[width=170mm]{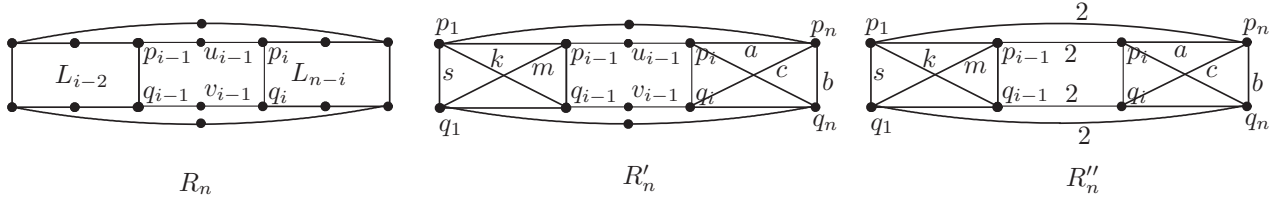}
  \caption{Graphs $R_n,\, R_n'$ and $R_n''$.}
\end{figure}

With the ordering of the vertices $(p_1,p_{i-1},p_i,p_n,q_1,q_{i-1},q_i,q_n)$, the discrete Laplacian matrix $L$ of the graph $R_n''$ is as follows:
$$
L=\left(
  \begin{array}{rrrrrrrr}
    K & -\frac{1}{m} & 0 & -\frac{1}{2} & -\frac{1}{s} & -\frac{1}{k} & 0 & 0 \\[2pt]
    -\frac{1}{m} & K & -\frac{1}{2} & 0 & -\frac{1}{k} & -\frac{1}{s} & 0 & 0 \\[2pt]
    0 & -\frac{1}{2} & S & -\frac{1}{a} & 0 & 0 & -\frac{1}{b} & -\frac{1}{c} \\[2pt]
    -\frac{1}{2} & 0 & -\frac{1}{a} & S & 0 & 0 & -\frac{1}{c} & -\frac{1}{b} \\[2pt]
    -\frac{1}{s} & -\frac{1}{k} & 0 & 0 & K & -\frac{1}{m} & 0 & -\frac{1}{2} \\[2pt]
    -\frac{1}{k} & -\frac{1}{s} & 0 & 0 & -\frac{1}{m} & k & -\frac{1}{2} & 0 \\[2pt]
    0 & 0 & -\frac{1}{b} & -\frac{1}{c} & 0 & -\frac{1}{2} & S & -\frac{1}{2} \\[2pt]
    0 & 0 & -\frac{1}{c} & -\frac{1}{b} & -\frac{1}{2} & 0 & -\frac{1}{a} & S
  \end{array}
 \right),
$$
where $K=\frac{1}{2}+\frac{1}{m}+\frac{1}{s}+\frac{1}{k}$ and $S=\frac{1}{2}+\frac{1}{a}+\frac{1}{b}+\frac{1}{c}$. Then we obtain the Moore-Penrose inverse $L^{+}$ of $L$ (see \cite{Bapat3}) as
$$ L^{+}=(L+\frac{1}{8}J)^{-1}-\frac{1}{8}J,$$
where $J$ is the matrix with all entries 1.

Next, together with Lemma \ref{lem2}, the result of $L^{+}$ and do some algebra by \cite{math} we obtain
\begin{align}
r(p_1,p_{i-1})=&\frac{1}{2}(\frac{1}{\frac{1}{4+\frac{ac}{a+c}}+\frac{1}{\frac{mk}{m+k}}}+\frac{1}{\frac{1}{4+\frac{ab}{a+b}}+\frac{1}{\frac{ms}{m+s}}}),\ \ \
r(p_1,q_{i-1})=\frac{1}{2}(\frac{1}{\frac{1}{4+\frac{ac}{a+c}}+\frac{1}{\frac{mk}{m+k}}}+\frac{1}{\frac{1}{4+\frac{bc}{b+c}}+\frac{1}{\frac{ks}{k+s}}}).\label{eq:5.8}
\end{align}

Note that (\ref{eq:5.5}), (\ref{eq:5.6}) and (\ref{eq:5.7}) describe the properties of the resistances along edges $w$, $u$ and $t$ in $K_4$. According to the structure of $R_n''$, we have
\begin{align*}
\frac{ac}{a+c}=C_{n-i},\ \frac{ab}{a+b}=B_{n-i},\ \frac{bc}{b+c}=A_{n-i},\
\frac{mk}{m+k}=C_{i-2},\ \frac{ms}{m+s}=B_{i-2},\ \frac{ks}{k+s}=A_{i-2}.
\end{align*}

Then we can change the subscript in (\ref{eq:5.8}) and simplify it as
$$r(p_1,p_i)=\frac{1}{2}(\frac{1}{\frac{1}{4+C_{n-i-1}}+\frac{1}{C_{i-1}}}+\frac{1}{\frac{1}{4+B_{n-i-1}}+\frac{1}{B_{i-1}}}),\ \ r(p_1,q_i)=\frac{1}{2}(\frac{1}{\frac{1}{4+C_{n-i-1}}+\frac{1}{C_{i-1}}}+\frac{1}{\frac{1}{4+A_{n-i-1}}+\frac{1}{A_{i-1}}}).$$

Do some algebra by \cite{math}, we have
\begin{align}
r(p_1,p_i)&=\frac{1+\alpha^{n}-\alpha^{n-i+1}-\alpha^{i-1}}{2\sqrt{2}(1-\alpha^{n})}+\frac{(n-i+1)(i-1)}{n},\label{eq:5.9}\\
r(p_1,q_i)&=\frac{1+\alpha^{n}+\alpha^{n-i+1}+\alpha^{i-1}}{2\sqrt{2}(1-\alpha^{n})}+\frac{(n-i+1)(i-1)}{n},\label{eq:5.10}
\end{align}
where $1\leq i\leq n$.

Next, our aim is to determine formulae for $r(u_1,u_i)$ and $r(u_1,v_i)$. In this case, let the 2-degree vertices in $R_n'$ absorb into $K_{4}$. Similar to the proof of Lemma {\ref{lem:5.1}}, we have the following result.
\begin{lem}\label{lem:5.2}
Assume that $L_{i-2}''=L_{i-2}+\{u_np_1,v_nq_1,u_{i-1}p_{i-1},v_{i-1}q_{i-1}\}$ , then we can transform $L_{i-2}''$ to $K_4$ where $x$, $y$ and $z$ are the resistances along edges in $K_4$ (see Fig. 13). This transformation yields equivalent circuits and the following hold.
\begin{align}
\frac{yz}{y+z}=2+A_{i-1}=-2\sqrt{2}+\frac{4\sqrt{2}}{1-\alpha^{i}},\ \
\frac{xy}{x+y}=2+B_{i-1}=-2\sqrt{2}+\frac{4\sqrt{2}}{1+\alpha^{i}},\ \
\frac{xz}{x+z}=2+C_{i-1}=2i.\label{eq:5.11}
\end{align}
\end{lem}
\begin{figure}[H]
\centering
\psfrag{y}{$u_n$}
\psfrag{z}{$p_1$}
\psfrag{6}{$v_1$}
\psfrag{7}{$q_1$}
\psfrag{r}{$u_{i-1}$}
\psfrag{9}{$p_{i-1}$}
\psfrag{8}{$q_{i-1}$}
\psfrag{m}{$v_{i-1}$}
\psfrag{2}{$L_{i-2}$}
\psfrag{h}{$\Rightarrow$}
\psfrag{x}{$u_n$}
\psfrag{5}{$v_n$}
\psfrag{u}{$x$}
\psfrag{w}{$y$}
\psfrag{t}{$z$}
\includegraphics[width=4in]{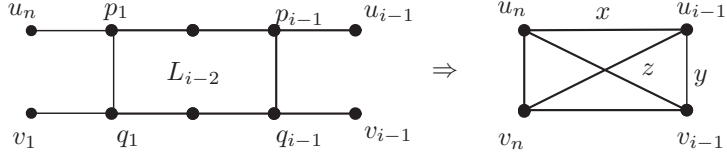}
\caption{Simplified circuit of $L_{i-2}''$.}
\end{figure}

Now we use Lemma \ref{lem:5.2} to simplify the circuit of $R_n$ as the form of $R_n'''$ (see Fig. 14).
\begin{figure}[H]
\centering
\psfrag{a}{$a$}
\psfrag{b}{$b$}
\psfrag{c}{$c$}
\psfrag{k}{$p_n$}
\psfrag{u}{$q_n$}
\psfrag{z}{$p_i$}
\psfrag{7}{$q_i$}
\psfrag{8}{$1$}
\psfrag{x}{$p_{i-1}$}
\psfrag{5}{$q_{i-1}$}
\psfrag{c}{$c$}
\psfrag{w}{$x$}
\psfrag{s}{$y$}
\psfrag{m}{$z$}
\psfrag{t}{$u_n$}
\psfrag{h}{$v_n$}
\psfrag{r}{$\Rightarrow$}
\psfrag{1}{$L_{i-2}$}
\psfrag{2}{$L_{n-i}$}
\psfrag{y}{$u_{i-1}$}
\psfrag{6}{$v_{i-1}$}
\psfrag{f}{$R_{n}$}
\psfrag{g}{$R_{n}'''$}
\psfrag{u}{$u_n$}
\psfrag{v}{$v_n$}
\psfrag{0}{$u_{i-1}$}
\psfrag{3}{$v_{i-1}$}
\psfrag{4}{$p_i$}
\psfrag{9}{$q_i$}
\includegraphics[width=6in]{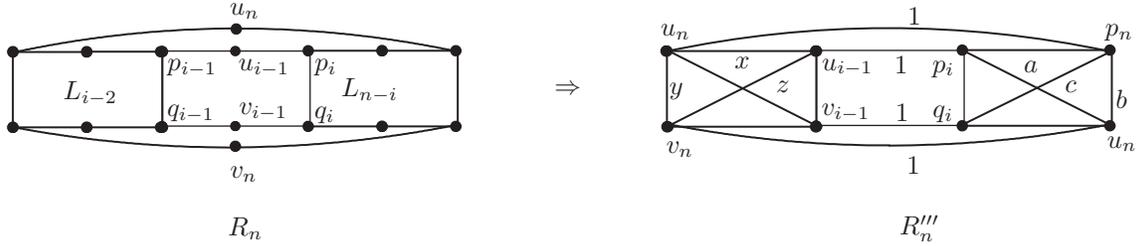}
\caption{$R_n\ \ \Rightarrow \ \ R_{n}'''$.}
\end{figure}
Similarly, we combine the Moore-Penrose inverse and Lemma {\ref{lem2}} to obtain that
\begin{align}
r_{R_n'''}(u_n,u_{i-1})=\frac{1}{2}(\frac{1}{\frac{1}{2+\frac{ac}{a+c}}+\frac{1}{\frac{xz}{x+z}}}+\frac{1}{\frac{1}{2+\frac{ab}{a+b}}+\frac{1}{\frac{xy}{x+y}}}),\ \
r_{R_n'''}(v_n,v_{i-1})=\frac{1}{2}(\frac{1}{\frac{1}{2+\frac{ac}{a+c}}+\frac{1}{\frac{xz}{x+z}}}+\frac{1}{\frac{1}{2+\frac{bc}{b+c}}+\frac{1}{\frac{yz}{y+z}}}).
\label{eq:5.12}
\end{align}
Then we use Eq. ({\ref{eq:5.11}}) to simplify Eq. ({\ref{eq:5.12}}) and get the formulae for $r_{R_n}(u_1,u_i)$ and $r_{R_n}(u_1,v_i)$ as
\begin{align}
r(u_1,u_i)=&r_{R_n'''}(u_n,u_{i-1})
=\frac{1}{2}(\frac{1}{\frac{1}{2+C_{n-i}}+\frac{1}{\frac{xz}{x+z}}}+\frac{1}{\frac{1}{2+B_{n-i}}+\frac{1}{\frac{xy}{x+y}}})\notag
\end{align}
\begin{align}
=&\frac{1+\alpha^{n}-\alpha^{n-i+1}-\alpha^{i-1}}{\sqrt{2}(1-\alpha^{n})}+\frac{(n-i+1)(i-1)}{n},\label{eq:5.13}\\
r(u_1,v_i)=&r_{R_n'''}(u_n,v_{i-1})
=\frac{1}{2}(\frac{1}{\frac{1}{2+C_{n-i}}+\frac{1}{\frac{xz}{x+z}}}+\frac{1}{\frac{1}{2+A_{n-i}}+\frac{1}{\frac{yz}{y+z}}})\notag\\
=&\frac{1+\alpha^{n}+\alpha^{n-i+1}+\alpha^{i-1}}{\sqrt{2}(1-\alpha^{n})}+\frac{(n-i+1)(i-1)}{n},\label{eq:5.14}
\end{align}
where $1\leq i\leq n$.

Finally, we try to find the effective resistances between $u_{i-1}$, $v_{i-1}$ and $p_1$. To solve the problem, we also use Lemma {\ref{lem:5.2}} to simplify $R_n$ as the form of $R_n''''$ (see Fig. 15). The parameter $a'$, $b'$ and $c'$ are similar to $x$, $y$ and $z$ in the previous case.
\begin{figure}[H]
\centering
\psfrag{a}{$a'$}
\psfrag{b}{$b'$}
\psfrag{c}{$c'$}
\psfrag{k}{$u_n$}
\psfrag{u}{$v_n$}
\psfrag{z}{$p_{i}$}
\psfrag{7}{$q_{i}$}
\psfrag{8}{$1$}
\psfrag{x}{$p_{i-1}$}
\psfrag{5}{$q_{i-1}$}
\psfrag{w}{$m$}
\psfrag{s}{$s$}
\psfrag{m}{$k$}
\psfrag{t}{$p_1$}
\psfrag{h}{$q_1$}
\psfrag{r}{$\Rightarrow$}
\psfrag{1}{$L_{i-2}$}
\psfrag{2}{$L_{n-i}$}
\psfrag{y}{$u_{i-1}$}
\psfrag{6}{$v_{i-1}$}
\psfrag{f}{$R_{n}$}
\psfrag{g}{$R_{n}''''$}
\psfrag{u}{$u_n$}
\psfrag{v}{$v_n$}
\psfrag{0}{$p_{i-1}$}
\psfrag{3}{$q_{i-1}$}
\psfrag{4}{$u_{i-1}$}
\psfrag{9}{$v_{i-1}$}
\includegraphics[width=6in]{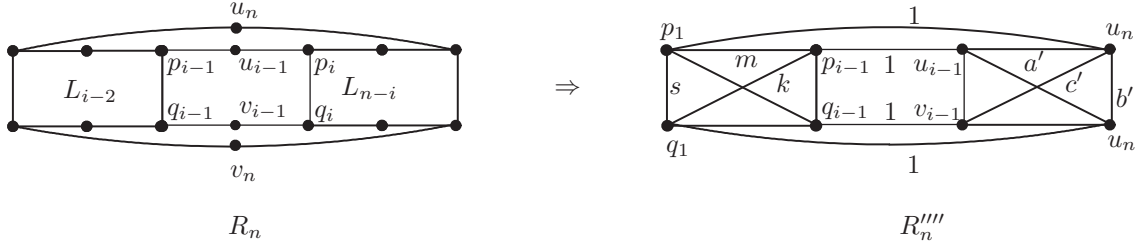}
\caption{$R_n\ \ \Rightarrow \ \ R_{n}''''$.}
\end{figure}

We use Lemmas {\ref{lem:5.1}} and {\ref{lem:5.2}} doing some equivalent deformation of parameters. Let
$$
\begin{array}{llll}
   &A'=\frac{a'c'}{a'+c'}=2(n-i+1), &B'=\frac{a'b'}{a'+b'}=-2\sqrt{2}+\frac{4\sqrt{2}}{1+\alpha^{n-i+1}}, &C'=\frac{b'c'}{b'+c'}=-2\sqrt{2}+\frac{4\sqrt{2}}{1-\alpha^{n-i+1}},  \\[5pt]
   &M'=\frac{mk}{m+k}=2(i-2), &S'=\frac{ms}{m+s}=-2-2\sqrt{2}+\frac{4\sqrt{2}}{1+\alpha^{i-1}}, &K'=\frac{ks}{k+s}=-2-2\sqrt{2}+\frac{4\sqrt{2}}{1-\alpha^{i-1}}.
\end{array}
$$
Then we obtain the formulae for $r_{R_n}(p_1,u_{i-1})$ and $r_{R_n}(p_1,v_{i-1})$ with respect to $A', B', C', M', S', K'$ as
\begin{align*}
r(p_1,u_{i-1})=&\frac{1}{2}\left[\frac{3+2A'+S'(2+A')+M'(2+S'+K')}{2+M'+A'}-\frac{(1+S')^{2}}{2+S'+B'}-\frac{1}{2+K'+C'}\right]\\
=&\frac{11+2\sqrt{2}-(20+14\sqrt{2})\alpha^{i}-(20-14\sqrt{2})\alpha^{n-i}-(5+2\sqrt{2})\alpha^{n}}{4\sqrt{2}(1-\alpha^{n})}-\frac{(2i-3)^{2}}{4n}+i-2-\sqrt{2},\\
r(p_1,v_{i-1})=&\frac{1}{2}\left[\frac{3+2A'+K'(2+A')+M'(2+K'+S')}{2+M'+A'}-\frac{1}{2+S'+B'}-\frac{(1+K')^{2}}{2+K'+C'}\right]\\
=&\frac{11+2\sqrt{2}+(20+14\sqrt{2})\alpha^{i}+(20-14\sqrt{2})\alpha^{n-i}-(5+2\sqrt{2})\alpha^{n}}{4\sqrt{2}(1-\alpha^{n})}-\frac{(2i-3)^{2}}{4n}+i-2-\sqrt{2},
\end{align*}
where $1\leq i-1\leq n-1$.

 That is
\begin{align}
r(p_1,u_{i})=&\frac{11+2\sqrt{2}-(20+14\sqrt{2})\alpha^{i+1}-(20-14\sqrt{2})\alpha^{n-i-1}-(5+2\sqrt{2})\alpha^{n}}{4\sqrt{2}(1-\alpha^{n})}-\frac{(2i-1)^{2}}{4n}+i-1-\sqrt{2},\label{eq:5.15}\\
r(p_1,v_{i})=&\frac{11+2\sqrt{2}+(20+14\sqrt{2})\alpha^{i+1}+(20-14\sqrt{2})\alpha^{n-i-1}-(5+2\sqrt{2})\alpha^{n}}{4\sqrt{2}(1-\alpha^{n})}-\frac{(2i-1)^{2}}{4n}+i-1-\sqrt{2},\label{eq:5.16}
\end{align}
where $1\leq i\leq n-1$.

Note that ({\ref{eq:5.15}}),({\ref{eq:5.16}}) hold under the condition $1\leq i< n$. However, these formulae are consistent with $i=n$ based on  $r(p_1,u_n)=r(p_1,u_1)$ and $r(p_1,v_n)=r(p_1,v_1)$.
 \subsection{\normalsize The maximum and minimum effective resistances in $R_n$}

In this subsection, we consider the extremal problems about the effective resistances in $R_n$. 
Note that $R_n$ is symmetric, hence for $r(p_j,p_i)$, $r(p_j,q_i)$, $r(u_j,u_i)$, $r(u_j,v_i)$, $r(p_j,u_i)$ and $r(p_j,v_i)$, we may fix $j=1$.
\begin{lem}\label{lem:6.1}
In the graph $R_n$, for any $1\leq i \leq n$, $r(p_1,p_i)$, $r(p_1,q_i)$, $r(u_1,u_i)$, $r(u_1,v_i)$, $r(p_1,u_i)$ and $r(p_1,v_i)$ are convex in i, i.e,
\begin{align*}
r&(p_1,p_2)<r(p_2,p_{3})<\dots <r(p_1,p_{\lfloor \frac{n+2}{2} \rfloor})
=r(p_1,p_{\lceil \frac{n+2}{2} \rceil})>\dots>r(p_1,p_n),\quad k> 0.\\
r&(p_1,q_1)<r(p_1,q_2)<\dots <r(p_1,q_{\lfloor \frac{n+2}{2} \rfloor})
=r(p_1,q_{\lceil \frac{n+2}{2} \rceil})>\dots>r(p_1,q_n),\quad k\geq 0.\\
r&(u_1,u_2)<r(u_2,u_{3})<\dots <r(u_1,u_{\lfloor \frac{n+2}{2} \rfloor})
=r(u_1,u_{\lceil \frac{n+2}{2} \rceil})>\dots>r(u_1,u_n),\quad k> 0.\\
r&(u_1,v_1)<r(u_1,v_2)<\dots <r(u_1,v_{\lfloor \frac{n+2}{2} \rfloor})
=r(u_1,v_{\lceil \frac{n+2}{2} \rceil})>\dots>r(u_1,v_n),\quad k\geq 0.\\
r&(p_1,u_1)<r(p_1,u_2)<\dots <r(p_1,u_{\lfloor \frac{n+1}{2} \rfloor})
=r(p_1,u_{\lceil \frac{n+1}{2} \rceil})>\dots>r(p_1,u_n),\quad k\geq 0.\\
r&(p_1,v_1)<r(p_1,v_2)<\dots <r(p_1,v_{\lfloor \frac{n+1}{2} \rfloor})
=r(p_1,v_{\lceil \frac{n+1}{2} \rceil})>\dots>r(p_1,v_n),\quad k\geq 0.
\end{align*}
\end{lem}
\begin{proof}
First, we try to proof this lemma for $r(p_1,p_i)$. That is to consider the sign of $r(p_1,p_{i+1})-r_{R_n}(p_1,p_i)$.
$$r(p_1,p_{i+1})-r(p_1,p_i)=\frac{(1-\alpha)\alpha^{-1-i}(\alpha^{2i}-\alpha^{n+1})}{2\sqrt{2}(1-\alpha^{n})}+\frac{n+1-2i}{n}.$$
By a direct calculation, we have $r(p_1,p_{i+1})  \geq  r(p_1,p_i)$ if $n\geq 2i-1$ and $r(p_1,p_{i+1})  <  r(p_1,p_i)$ otherwise.

Now we consider $r(p_1,u_i)$. In fact,
$$ r_{R_n}(p_1,u_{i+1})-r_{R_n}(p_1,u_i)=\frac{\alpha^{-i}(\alpha^{2i}-n)}{1-\alpha^{n}}+\frac{n-2i}{n}.$$
Hence,
Then we claim that $r(p_1,u_{i+1})  \geq  r(p_1,u_i)$ if $n\geq 2i$ and $r(p_1,u_{i+1})  <  r(p_1,u_i)$ otherwise.

Similarly, we may show the following cases, which are omitted here.
\end{proof}

Using Lemma {\ref{lem:6.1}}, we can obtain the following result.
\begin{thm}
In the graph $R_n$, for fixed $n\geq3$ and any vertex $a,b\in V_{R_n}$,
$r(u_1,v_{\lfloor \frac{n+2}{2} \rfloor})\geq r(a,b)\geq r(p_1,q_1).$
\end{thm}
\begin{proof}
We consider firstly the maximum effective resistances. By Lemma {\ref{lem:6.1}}, it suffices to determine
$$\max\{r(p_1,p_{\lfloor \frac{n+2}{2} \rfloor}), r(p_1,q_{\lfloor \frac{n+2}{2} \rfloor}), r(u_1,u_{\lfloor \frac{n+2}{2} \rfloor}), r(u_1,v_{\lfloor \frac{n+2}{2} \rfloor}), r(p_1,u_{\lfloor \frac{n+1}{2} \rfloor}), r(p_1,v_{\lfloor \frac{n+1}{2} \rfloor})\}.$$

If $n$ is even, then
\begin{align*}
r(p_1,p_{\lfloor \frac{n+2}{2} \rfloor})=&\frac{1}{2\sqrt{2}}(\frac{2}{1+\alpha^{n/2}}+n-1),\ \ \ \ \ \
r(p_1,q_{\lfloor \frac{n+2}{2} \rfloor})=\frac{1}{2\sqrt{2}}(\frac{2}{1-\alpha^{n/2}}+n-1),\\
r(u_1,u_{\lfloor \frac{n+2}{2} \rfloor})=&\frac{1-\alpha^{n/2}}{\sqrt{2}(1+\alpha^{n/2})}+\frac{n}{4},\ \ \ \ \ \ \ \ \ \ \ \ \,
r(u_1,v_{\lfloor \frac{n+2}{2} \rfloor})=\frac{1+\alpha^{n/2}}{\sqrt{2}(1-\alpha^{n/2})}+\frac{n}{4},\\
r(p_1,u_{\lfloor \frac{n+1}{2} \rfloor})=&\frac{3-8\alpha^{n/2}+3\alpha^{n}}{4\sqrt{2}(1-\alpha^{n})}-\frac{1}{4n}+\frac{n}{4},\
r(p_1,v_{\lfloor \frac{n+1}{2} \rfloor})=\frac{3+8\alpha^{n/2}+3\alpha^{n}}{4\sqrt{2}(1-\alpha^{n})}-\frac{1}{4n}+\frac{n}{4}.
\end{align*}
This gives
\begin{align*}
r(u_1,v_{\lfloor \frac{n+2}{2} \rfloor})-r(p_1,p_{\lfloor \frac{n+2}{2} \rfloor})&=
\frac{1+6\alpha^{n/2}+\alpha^{n}}{2\sqrt{2}(1-\alpha^{n})}>0,\\
r(u_1,v_{\lfloor \frac{n+2}{2} \rfloor})-r(p_1,q_{\lfloor \frac{n+2}{2} \rfloor})&=
\frac{1+\alpha^{n/2}}{2\sqrt{2}(1-\alpha^{n/2})}>0,\\
r(u_1,v_{\lfloor \frac{n+2}{2} \rfloor})-r(u_1,u_{\lfloor \frac{n+2}{2} \rfloor})&=
\frac{2\sqrt{2}\alpha^{n/2}}{1-\alpha^{n}}>0,\\
r(u_1,v_{\lfloor \frac{n+2}{2} \rfloor})-r(p_1,u_{\lfloor \frac{n+1}{2} \rfloor})&=
\frac{-4-7\sqrt{2}+16\alpha^{n/2}+(4+9\sqrt{2})\alpha^{n}}{8(1-\alpha^{n})}+\frac{1}{2}+\sqrt{2}-\frac{1}{4n}>0,\\
r(u_1,v_{\lfloor \frac{n+2}{2} \rfloor})-r(p_1,v_{\lfloor \frac{n+1}{2} \rfloor})&=
\frac{-4-7\sqrt{2}+(4+9\sqrt{2})\alpha^{n}}{8(1-\alpha^{n})}+\frac{1}{2}+\sqrt{2}-\frac{1}{4n}>0.
\end{align*}
Thus, $r(u_1,v_{\lfloor \frac{n+2}{2} \rfloor})$ is the maximum resistance for even $n$. Similarly, $r(u_1,v_{\lfloor \frac{n+2}{2} \rfloor})$ is the maximum resistance for odd $n$. We omit the procedure here. 

Now we consider the minimum effective resistance. By Lemma {\ref{lem:6.1}}, it suffices to determine
$$
\min\{r(p_1,p_2),\, r(p_1,q_1),\, r(u_1,v_1),\, r(u_1,u_2),\, r(p_1,u_1),\, r(p_1,v_1)\}.
$$

In fact, by a direct calculation one has
$$
\begin{array}{lll}
  r(p_1,p_2)=\frac{4-\sqrt{2}-(4+\sqrt{2})\alpha^{n}}{2(1-\alpha^{n})}-\frac{1}{n}, & r(p_1,q_1)=\frac{1+\alpha^{n}}{\sqrt{2}(1-\alpha^{n})}, & r(u_1,u_2)=\frac{2-\sqrt{2}-(2+\sqrt{2})\alpha^{n}}{1-\alpha^{n}}+\frac{n-1}{n}, \\
   r(u_1,v_1)=\frac{\sqrt{2}(1+\alpha^{n})}{1-\alpha^{n}}, & r(p_1,u_1)=\frac{4\sqrt{2}-1-(1+4\sqrt{2})\alpha^{n}}{4\sqrt{2}(1-\alpha^{n})}-\frac{1}{4n},  &  r(p_1,v_1)=\frac{7(1+\alpha^{n})}{4\sqrt{2}(1-\alpha^{n})}-\frac{1}{4n}.
\end{array}
$$
By a direct calculation, one has $
\min\{r(p_1,p_2),\, r(p_1,q_1),\, r(u_1,v_1),\, r(u_1,u_2),\, r(p_1,u_1),\, r(p_1,v_1)\}=r(p_1,q_1),
$
as desired.
\end{proof}

At the end of this section, we turn to the asymptotic properties of resistance distances in $R_n$.
\begin{thm}
\begin{wst}
\item[{\rm (i)}] For all fixed $i$ and $j$, one has
$$\lim_{n\to \infty}r(p_i,p_j)=i-1+\frac{1-\alpha^{i-1}}{2\sqrt{2}},\quad \quad
  \lim_{n\to \infty}r(q_i,p_j)=i-1+\frac{1+\alpha^{i-1}}{2\sqrt{2}}.$$
\item[{\rm (ii)}] $\lim\limits_{n \to \infty}\frac{1}{n}r(p_1,p_{\lfloor \frac{n+2}{2} \rfloor})=\lim\limits_{n \to \infty}\frac{1}{n}r(p_1,q_{\lfloor \frac{n+2}{2} \rfloor})=\frac{1}{4}.$
\item[{\rm (iii)}]
  $\lim\limits_{n \to \infty} [r_{R_{n+1}}(p_1,p_{\lfloor \frac{n+3}{2} \rfloor})-r_{R_n}(p_1,p_{\lfloor \frac{n+2}{2} \rfloor})]=
  \lim\limits_{n \to \infty} [r_{R_{n+1}}(p_1,q_{\lfloor \frac{n+3}{2} \rfloor})-r_{R_n}(p_1,q_{\lfloor \frac{n+2}{2} \rfloor})]=\frac{1}{4}.$
\item[{\rm (iv)}]
  $\lim\limits_{n \to \infty}r(p_1,p_2)=2-\frac{\sqrt{2}}{2},
  \quad \lim\limits_{n \to \infty}r(p_1,q_1)=\frac{\sqrt{2}}{2}.$
\end{wst}
\end{thm}
\begin{proof}
(i)-(iv) follow directly from ({\ref{eq:5.9}}) and (\ref{eq:5.10}).
\end{proof}
By a similar discussion, it is not difficult to determine the limit value on the resistance distance between any other pair of vertices in $R_n$ as $n\rightarrow \infty$. Here we omit the contents.
\section{\normalsize The Kirchhoff indices of $L_n$ and $R_n$}\setcounter{equation}{0}

In this section, we determine the formulae for the Kirchhoff indices of $L_n$ and $R_n$. Recall that \textit{Kirchhoff index} of a graph $G$, $K\!f(G)$ is defined \cite{Klein} as follows:
\begin{align}\label{eq33}
 K\!f(G)=\sum_{\{v,u\}\subseteq V_{G}}r(u,v).
\end{align}
\begin{thm}
Let $L_n$ be a linear hexagonal chain. Then
$$K\!f(L_n)=\frac{(1+2n)(21-6\sqrt{2}+2(8+9\sqrt{2})n-\alpha^{2n+2}(21+6\sqrt{2}+2(8-9\sqrt{2})n))}{12(1-\alpha^{2n+2})}+\frac{4}{3}(1+2n)n^2.$$
\end{thm}
\begin{proof}
According to the definition of Kirchhoff index, together with ({\ref{eq:3.9}})-({\ref{eq:3.14}}), we have
\begin{align}
K\!f(L_n)=&\frac{1}{2}\sum_{a, b\in V(L_n)} r(a,b)\notag\\
=&\sum_{i=0}^{n}r(p_i,q_i)+\sum_{j=0}^{n}\sum_{i=j+1}^{n}[r(p_i,p_j)+r(q_i,q_j)+r_(p_i,q_j)+r(q_i,p_j)]+\sum_{i=0}^{n-1}r(u_i,v_i)+\sum_{j=0}^{n-1}\sum_{i=j+1}^{n-1}\left[r(u_i,u_j)\right.\notag\\
&\left.+r(v_i,v_j)+r(u_i,v_j)+r(v_i,u_j)\right]+2\sum_{j=0}^{n-1}\sum_{i=j}^{n-1}[r(u_i,p_j)+r(v_i,q_j)+r(v_i,p_j)+r(u_i,q_j)]\notag\\
=&\sum_{i=0}^{n}r(p_i,q_i)+2\sum_{j=0}^{n}\sum_{i=j+1}^{n}[r(p_i,p_j)+r_(p_i,q_j)]+\sum_{i=0}^{n-1}r(u_i,v_i)+2\sum_{j=0}^{n-1}\sum_{i=j+1}^{n-1}[r(u_i,u_j)+r(u_i,v_j)]\notag\\
&+4\sum_{j=0}^{n-1}\sum_{i=j}^{n-1}[r(u_i,p_j)+r(v_i,p_j)]\notag\\
=&\frac{(1+2n)(21-6\sqrt{2}+2(8+9\sqrt{2})n-\alpha^{2n+2}(21+6\sqrt{2}+2(8-9\sqrt{2})n))}{12(1-\alpha^{2n+2})}+\frac{4}{3}(1+2n)n^2,\label{eq:7.1}
\end{align}
where (\ref{eq:7.1}) is obtained by doing some algebra through \cite{math}.
\end{proof}

Similarly, we use \cite{math} to obtain the following result by ({\ref{eq:5.9}}), ({\ref{eq:5.10}}), ({\ref{eq:5.13}})-({\ref{eq:5.16}}).
\begin{thm}
Let $R_n$ be a hexagonal cylinder chain. Then
\begin{align*}
K\!f(R_n)=&
\frac{4n^{3}-n}{3}+3\sqrt{2} n^{2}\frac{1+\alpha^{n}}{1-\alpha^{n}}.
\end{align*}
\end{thm}

In view of Theorems 5.1 and 5.2, the next corollary follows directly.
\begin{cor}
Let $L_n$ be a linear hexagonal chain and $R_n$ be a hexagonal cylinder chain. Then
$$\lim_{n\to +\infty}\frac{K\!f(R_n)}{K\!f(L_n)}=\frac{1}{2}.$$
\end{cor}

\end{document}